\documentclass{amsart}

\numberwithin{equation}{section}

\usepackage{amssymb}
\usepackage{enumerate, xspace}
\usepackage[matrix,arrow]{xy}
\usepackage{comment}

\newtheorem{thm}{Theorem}[section]
\newtheorem{lem}[thm]{Lemma}
\newtheorem{cor}[thm]{Corollary}

\newtheorem{prop}[thm]{Proposition}

\newtheorem{defin}[thm]{Definition}

\newtheorem{conv}[thm]{Convention}

\newtheorem{rem}[thm]{Remark}
\newtheorem{dft}[thm]{Draft Remark}

\newcommand\B{{\mathcal B}}
\newcommand\Co{{\mathcal C}}
\newcommand\D{{\mathcal D}}
\newcommand\E{{\mathcal E}}
\newcommand\F{{\mathcal F}}
\newcommand\G{{\mathcal G}}

\newcommand\Lp{{\mathcal L}}
\newcommand\Or{{\mathcal O}}
\newcommand\M{{\mathcal M}}
\newcommand\Se{{\mathcal S}}
\newcommand\T{{\mathcal T}}

\newcommand\A{{\mathbb A}}

\newcommand\C{{\mathbb C}}

\newcommand\N{{\mathbb N}}

\newcommand\R{{\mathbb R}}

\newcommand\V{{\mathbb V}}
\newcommand\Z{{\mathbb Z}}

\newcommand\ve{\varepsilon}

\newcommand\vf{\varphi}

\newcommand{\ra}{\varrho}

\newcommand{\coe}{\iota}
\newcommand{\vectfield}{\mathcal{V}}
\newcommand{\trunc}{\pi}
\newcommand{\Weight}{\mathcal{W}}
 
\newcommand{\norm}[1]{\left\| #1 \right\|}
\newcommand{\chart}{\psi}
\newcommand{\expansion}{\kappa}
\newcommand{\dd}{\, {\rm d}}
\DeclareMathOperator{\pr}{pr} \DeclareMathOperator{\diam}{diam}
\newcommand{\FF}{\mathcal{F}}

\newcommand{\Ptop}{P_{\rm top}(\bar \phi)}
\newcommand{\CC}{\mathcal{C}}
\newcommand{\ora}{\bar\ra}
\newcommand{\st}{\;:\;}
\DeclareMathOperator{\dist}{dist}
\newcommand{\betriem}{\alpha_{r}}
\newcommand{\betlim}{\alpha_0}
\DeclareMathOperator{\dLeb}{dLeb}
\excludecomment{dft}

\begin{document}

\title[Hyperbolic sets]{Compact locally maximal hyperbolic sets for
smooth maps: fine statistical properties}
\author{S\'{e}bastien Gou\"{e}zel and Carlangelo Liverani}
\address{S\'{e}bastien Gou\"{e}zel\\
IRMAR\\
Universit\'{e} de Rennes 1\\
Campus de Beaulieu, b\^{a}timent 22\\
35042 Rennes Cedex, France.}
\email{sebastien.gouezel@univ-rennes1.fr}
\address{Carlangelo Liverani\\
Dipartimento di Matematica\\
II Universit\`{a} di Roma (Tor Vergata)\\
Via della Ricerca Scientifica, 00133 Roma, Italy.}
\email{liverani@mat.uniroma2.it}
\date{June 28, 2006}
\begin{abstract}
Compact locally maximal hyperbolic sets are studied via
geometrically defined functional spaces that take advantage of the
smoothness of the map in a neighborhood of the hyperbolic set. This
provides a self-contained theory that not only reproduces all the
known classical results but  gives also new insights on the
statistical properties of these systems.
\end{abstract}
\thanks{We wish to thank A.Avila, V.Baladi, D.Dolgopyat, D.Ruelle  and
M.Tsujii for helpful discussions. We acknowledge the support of the
Institut Henri Poincar\'{e} where this work was started (during the
trimester {\sl Time at Work}), the GDRE Grefi-Mefi and the M.I.U.R.
(Cofin 05-06 PRIN 2004028108) for partial support.} \maketitle

%%%PAPER
\section{Introduction}\label{sec:zero}

The ergodic properties of uniformly hyperbolic maps can be described
as follows. If $T$ is a topologically mixing map on a compact
locally maximal hyperbolic set $\Lambda$ belonging to some smooth
manifold $X$, and $\bar\phi : \Lambda \to \R$ is a H\"{o}lder continuous
function, then there exists a unique probability measure
$\mu_{\bar\phi}$ maximizing the variational principle with respect
to $\bar\phi$ (the {\sl Gibbs measure with potential $\bar\phi$}).
Moreover, this measure enjoys strong statistical properties
(exponential decay of correlations, central and local limit
theorem...). When $\Lambda$ is an attractor and the potential
$\bar\phi$ is the jacobian of the map in the unstable direction (or,
more generally, a function which is cohomologous to this one), then
the measure $\mu_{\bar\phi}$ is the so-called SRB measure, which
describes the asymptotic behavior of Lebesgue-almost every point in
a neighborhood of $\Lambda$.

The proof of these results, due among others to Anosov, Margulis,
Sinai, Ruelle, Bowen, is one of the main accomplishments of the
theory of dynamical systems in the 70's. The main argument of their
proof is to \emph{code} the system, that is, to prove that it is
semiconjugate to a subshift of finite type, to show the
corresponding results for subshifts (first unilateral, and then
bilateral), and to finally go back to the original system. These
arguments culminate in Bowen's monograph \cite{bowen}, where all the
previous results are proved. Let us also mention another approach,
using \emph{specification}, which gives existence and uniqueness of
Gibbs measures (but without exponential decay of correlations or
limit theorems) through purely topological arguments
\cite{bowen:specification}.

These methods and results have proved very fruitful for a manifold
of problems. However, problems and questions of a new type have
recently emerged, such as
  \begin{itemize}
  \item Strong statistical stability w.r.t. smooth or random perturbations;
  \item Precise description of the correlations;
  \item Relationships between dynamical properties and the zeroes of the zeta function in a large
  disk.
  \end{itemize}
It is possible to give partial answers to these questions using
coding (see e.g.~\cite{ ruelle_dist, haydn:zeta,
ruelle:differentation, pollicott:stability}), but their range is
limited by the H\"{o}lder continuity of the foliation: the coding map
can be at best H\"{o}lder continuous, and necessarily loses information
on the smoothness properties of the transformation.

Recently, \cite{bkl:spectre_anosov} introduced a more geometric
method to deal with these problems, in the case of the SRB measure.
It was still limited by the smoothness of the foliation, but it
paved the way to further progress. Indeed, Gou\"{e}zel and Liverani
could get close to optimal answers to the first two questions (for
the SRB measure of an Anosov map) in \cite{gouezel_liverani}. Baladi
and Tsujii finally reached the optimal results in
\cite{baladi:Cinfty, bt:aniso} for these two questions (for the SRB
measure of an hyperbolic attractor). A partial answer to the last
question was first given in \cite{liverani:zeta} and a complete
solution will appear in the paper \cite{bt:zeta}. See also the paper
\cite{liverani-tsujii:zeta} for a very simple, although non optimal,
argument.

The technical approach of these papers is as follows: they introduce
spaces $\B$ of distributions, and an operator $\Lp : \B \to \B$ with
good spectral properties such that, for all smooth functions
$\psi_1,\psi_2$ and all $n\in \N$,
  \begin{equation}
  \label{def_Leb}
  \int \psi_1 \cdot \psi_2\circ T^n \dLeb = \langle \Lp^n( \psi_1
  \dLeb), \psi_2 \rangle.
  \end{equation}
The operator $\Lp$ has a unique fixed point, which corresponds to
the SRB measure of the map $T$. The correlations are then given by
the remaining spectral data of $\Lp$. In addition, abstract spectral
theoretic arguments imply precise results on perturbations of $T$ or
zeta functions.

\medskip

In this paper, we extend to the setting of Gibbs measures the
results of \cite{gouezel_liverani}. This extension is not
straightforward for the following reasons. First, the previous
approaches for the SRB measure rely on the fact that there is
already a reference measure to work with, the Lebesgue measure. For
a general (yet to be constructed) Gibbs measure, there is no natural
analogous of \eqref{def_Leb} which could be used to define the
transfer operator $\Lp$. The technical consequence of this fact is
that our space will not be a space of distribution on the whole
space, rather a family of distributions on stable (or close to
stable) leaves. Second, the SRB measure corresponds to a potential
$\bar\phi_u$ -- minus the logarithm of the unstable jacobian, with
respect to some riemannian metric -- which is in general \emph{not}
smooth, while we want our spaces to deal with very smooth objects.
Notice however that $\bar\phi_u$ is cohomologous to a function which
can be written as $\phi(x,E^s(x))$ where $\phi$ is a smooth function
on the grassmannian of $d_s$ dimensional subspaces of the tangent
bundle $\T X$.\footnote{ Take $\phi(x,E)=\log( \det
DT(x)_{\upharpoonright_E})- \log( \det DT(x))$, where $\det$
indicates the jacobian with respect to the given riemannian metric.
Let $\bar \phi(x)= \phi(x, E^s(x))$ for $x\in\Lambda$. Since the
angle between the stable and unstable direction is bounded from
below, $\sum_{k=0}^{n-1} (\bar \phi\circ T^k - \bar \phi_u\circ
T^k)$ is uniformly bounded on $\Lambda$. By Livsic theorem, this
implies that $\bar\phi$ is cohomologous to $\bar\phi_u$. In
particular, they give rise to the same Gibbs measure.}
 This is the kind
of potential we will deal with.

The elements of our Banach space $\B$ will thus be objects ``which
can be integrated along small submanifolds of dimension $d_s$''
(where $d_s$ is the dimension of the stable manifolds). The first
idea would be to take for $\B$ a space of differential forms of
degree $d_s$. However, if $\alpha$ is such a form and $\phi$ is a
potential as above, then $e^{\phi} \alpha$ is not a differential
form any more. Hence, we will have to work with more general
objects. Essentially, the elements of $\B$ are objects which
associate, to any subspace $E$ of dimension $d_s$ of the tangent
space, a volume form on $E$. Such an object can be integrated along
$d_s$ dimensional submanifolds, as required, and can be multiplied
by $e^{\phi}$. We define then an operator $\Lp$ on $\B$ by $\Lp
\alpha = T_*( e^{\phi} \pi \alpha)$ where $\pi$ is a truncation
function (necessary to keep all the functions supported in a
neighborhood of $\Lambda$, if $\Lambda$ is not an attractor), and
$T_*$ denotes the (naturally defined) push-forward of an element of
$\B$ under $T$. We will construct on $\B$ norms for which $\Lp$ has
a good spectral behavior, in Section \ref{sec:one}.

The main steps of our analysis are then the following.
\begin{enumerate}
\item Prove a Lasota-Yorke inequality for $\Lp$ acting on $\B$, in
Lemma \ref{lem:LY} (by using the preliminary result in Lemma
\ref{MainDynamicalInequality}). This implies a good spectral
description of $\Lp$ on $\B$: the spectral radius is some abstract
quantity $\ra$, yet to be identified, and the essential spectral
radius is at most $\sigma \ra$ for some small constant $\sigma$,
related to the smoothness of the map. See Proposition
\ref{thm:SpectralGap} and Corollary \ref{cor:lowerbound}.
\item
In this general setting, we analyze superficially the peripheral
spectrum (that is, the eigenvalues of modulus $\ra$), in Subsection
\ref{sec:PerSpectrum}. We prove that $\ra$ is an eigenvalue, and
that there is a corresponding eigenfunction $\betlim$ which induces
a measure on $d_s$ dimensional submanifolds (Lemma
\ref{lem:AreMeasures}). This does not exclude the possibility of
Jordan blocks or strange eigenfunctions.
\item
In the topologically mixing case, we check that $\betlim$ is fully
supported. By some kind of bootstrapping argument, this implies that
$\norm{\Lp^n} \leq C \ra^n$, i.e., there is no Jordan block.
Moreover, there is no other eigenvalue of modulus $\ra$ (Theorem
\ref{thm:PeripheralSpectrum}).
\item
The adjoint of $\Lp$, acting on $\B'$, has an eigenfunction $\ell_0$
for the eigenvalue $\ra$. The linear form $\vf \mapsto \ell_0( \vf
\betlim)$ is in fact a measure $\mu$, this will be the desired Gibbs
measure. Moreover, the correlations of $\mu$ are described by the
spectral data of $\Lp$ acting on $\B$, as explained in Section
\ref{sec:description}.
\item
Finally, in Section \ref{sec:variational} we prove that the
dynamical balls have a very well controlled measure (bounded from
below and above), see Proposition \ref{lem:MeasureDynBalls}. This
yields $\ra=\Ptop$ and the fact that $\mu$ is the unique equilibrium
measure (Theorem \ref{thm:IsGibbs}).
\end{enumerate}

It is an interesting issue to know whether there can indeed be
Jordan blocks in the non topologically transitive case (this is not
excluded by our results). The most interesting parts of the proof
are probably the Lasota-Yorke estimate and the exclusion of Jordan
blocks. Although the core of the argument is rather short and
follows very closely the above scheme, the necessary presence of the
truncation function induces several technical complications, which
must be carefully taken care of and cloud a bit the overall logic.
Therefore, the reader is advised to use the previous sketch of proof
to find her way through the rigorous arguments. Note that the paper
is almost completely self-contained, it only uses the existence and
continuity of the stable and unstable foliation (and not their
H\"{o}lder continuity nor their absolute continuity).

In addition, note that the present setting allows very precise
answers to the first of the questions posed at the beginning of this
introduction thanks to the possibility of applying the perturbation
theory developed in \cite[section 8]{gouezel_liverani} and based on
\cite{keller_liverani}. Always in the spirit to help the reader we
will give a flavor of such possibilities in Section
\ref{sec:applications} together with some obvious and less obvious
examples to which our theory can be applied. 
In particular, in Proposition 8.1 we provide nice formulae for the 
derivative of the
topological pressure and the Gibbs measure in the case of systems 
depending smoothly
on a parameter.\footnote{Note that the formulae are in terms 
of exponentially
converging sums, hence they can be easily used to actually compute the above
quantities within a given precision.}
Finally, a technical
section (Section \ref{sec:LeafMeasures}) on the properties of
conformal leafwise measures is added both for completeness and
because of its possible interest as a separate result.

\begin{rem}
Let us point out that, although we follow the strategy of
\cite{gouezel_liverani}, similar results can be obtained also by
generalizing the Banach spaces in \cite{bt:aniso} (M.Tsujii, private
communication).
\end{rem}

\medskip

To conclude the introduction let us give the description we obtain
for the correlation functions. We consider an open set $U\subset X$
and a map  $T\in\Co^r(U,X)$,\footnote{\label{foot:Cr}Here, and in
the following, by $\Co^r$ we mean the Banach space of functions
continuously differentiable $\lfloor r\rfloor$ times, and with the
$\lfloor r\rfloor$th derivative H\"{o}lder continuous of exponent
$r-\lfloor r\rfloor$. Such a space is equipped with a norm
$|\cdot|_{\Co^r}$ such that $|fg|_{\Co^r}\leq
|f|_{\Co^r}|g|_{\Co^r}$, that is $(\Co^r,|\cdot|_{\Co^r})$ is a
Banach algebra. For example, if $r\in\N$, $|f|_{\Co^r}:=\sup_{k\leq
r}|f^{(k)}|_\infty2^{r-k}$ will do.} diffeomorphic on its image (for
some real $r>1$). Suppose further that $\Lambda:=\bigcap_{n\in\Z}T^n
U$ is non empty and compact. Finally, assume that $\Lambda$ is a
hyperbolic set for $T$. Such a set is a \emph{compact locally
maximal hyperbolic set}. Let $\lambda>1$ and $\nu<1$ be two
constants, respectively smaller than the minimal expansion of $T$ in
the unstable direction, and larger than the minimal contraction of
$T$ in the stable direction.

Denote by $\Weight^0$ the set of $\Co^{r-1}$ function $\phi$
associating, to each $x\in U$ and each $d_s$ dimensional subspace of
the tangent space $\T_x X$ at $x$, an element of $\R$. Denote by
$\Weight^1$ the set of $\Co^r$ functions $\phi:U\to \R$. For $x\in
\Lambda$, set $\bar\phi(x)=\phi(x,E^s(x))$ in the first case, and
$\bar\phi(x)=\phi(x)$ in the second case. This is a H\"{o}lder
continuous function on $\Lambda$. Assume that the restriction of $T$
to $\Lambda$ is topologically mixing.

\begin{thm}
\label{DescribesCorrelations} Let $\phi \in \Weight^{\coe}$ for some
$\coe \in \{0,1\}$. Let $p\in \N^*$ and $q\in\R_+^*$ satisfy
$p+q\leq r-1+\coe$ and $q\geq \coe$. Let $\sigma>
\max(\lambda^{-p},\nu^q)$.\footnote{\label{pressure}In fact, one can
obtain better bounds by considering $T^n$, for large $n$, instead of
$T$. We will not indulge on such subtleties to keep the exposition
as simple as possible.} Then there exists a unique measure $\mu$
maximizing the variational principle for the potential
$\bar\phi$,\footnote{Of course, this is nothing else than the
classical Gibbs measure associated to the potential $\bar\phi$.} and
there exist a constant $C>0$, a finite dimensional space $F$, a
linear map $M:F\to F$ having a simple eigenvalue at $1$ and no other
eigenvalue with modulus $\geq 1$, and two continuous mappings
$\tau_1: \Co^{p}(U) \to F$ and $\tau_2: \Co^{q}(U) \to F'$ such
that, for all $\psi_1\in \Co^{p}(U)$, $\psi_2 \in \Co^q(U)$ and for
all $n\in \N$,
  \begin{equation}
  \label{eq:Correlations}
  \left| \int \psi_1 \cdot \psi_2\circ T^n \dd\mu- \tau_2(\psi_2) M^n
  \tau_1(\psi_1)\right| \leq C \sigma^n |\psi_2|_{\Co^{q}(U)}
  |\psi_1|_{\Co^{p}(U)}.
  \end{equation}
\end{thm}
The coefficients of the maps $\tau_1$ and $\tau_2$ are therefore
distributions of order at most $p$ and $q$ respectively, describing
the decay of correlations of the functions.  They extend the Gibbs
distributions of \cite{ruelle_dist} to a higher smoothness setting.

When $T$ is $\Co^\infty$, we can take $p$ and $q$ arbitrarily large,
and get a description of the correlations up to an arbitrarily small
exponential error term. The SRB measure corresponds to a potential
in $\Weight^0$, as explained above, and the restriction on $p,q$ is
$p+q \leq r-1$, which corresponds to the classical Kitaev bound
\cite{kitaev:fredholm}.\footnote{In some cases, our bound is not
optimal since $p$ is restricted to be an integer.} Surprisingly,
when the weight function belongs to $\Weight^1$, we can get up to
$p+q=r$. In some sense, the results are better for maximal entropy
measures than for SRB measures!

It is enlightening to consider our spaces for expanding maps, that
is, when $d_s=0$. In this case, ``objects that can be integrated
along stable manifolds'' are simply objects assigning a value to a
point, i.e., functions. Our Banach space $\B^{p,q}$ becomes the
space of usual $\Co^p$ functions, and we are led to the results of
Ruelle in \cite{ruelle:fredholm}.

\tableofcontents

\section{The functional spaces}\label{sec:one}
Consider a $\Co^{r}$ differentiable manifold $X$. We start with few
preliminaries.

\subsection{A touch of functional analysis}\label{subsec:space}
To construct the functional spaces we are interested in, we will use
an abstract construction that applies to each pair $\V,\Omega$,
where $\V$ is a complex vector space and $\Omega\subset \V'$ is a
subset of the (algebraic) dual with the property
$\sup_{\ell\in\Omega}|\ell(h)|<\infty$ for each $h\in\V$. In such a
setting we can define a seminorm on $\V$ by
  \begin{equation}
  \norm{h}:=\sup_{\ell\in\Omega}|\ell(h)|.
  \end{equation}
If we call $\B$ the completion of $\V$ with respect to $\|\cdot\|$,
we obtain a Banach space. Note that, by construction, $\Omega$
belongs to the unit ball of $\B'$. When $\|\cdot\|$ is a norm on
$\V$, i.e., $\V_0:= \bigcap_{\ell \in \Omega} \ker \ell$ is reduced
to $\{0\}$, then $\V$ can be identified as a subspace of $\B$. In
general, however, there is only an inclusion of the algebraic
quotient $\V / \V_0$ in $\B$.

\begin{comment}
More can be said. For example, the following holds.
\begin{lem}
\label{lem:dualball} For each $\ell\in\B'$ there exits a sequence
$\ell_n$ in the span of $\Omega$ that converges to $\ell$ in the
weak-* topology.
\end{lem}
The proof of the above lemma, together with some further discussion
on the dual can be found in Section \ref{sec:dual}.
\end{comment}

\subsection{Differential geometry beyond forms}
Let $\G$ be the Grassmannian of $d_s$ dimensional oriented subspaces
of the tangent bundle $\T X$ to $X$. On it we can construct the
complex line bundle $\E:=\{(x,E,\omega)\st (x,E)\in\G,
\omega\in\bigwedge^{d_s}E'\otimes \C\}$. We can then consider the
vector space $\Se$ of the $\Co^{r-1}$ sections of the line bundle
$\E$. The point is that for each $\alpha\in\Se$, each $d_s$
dimensional oriented $\Co^1$ manifold $W$ and each
$\vf\in\Co^0(W,\C)$, we can define an integration of $\vf$ over $W$
as if $\alpha$ was a usual differential form, by the formula
  \begin{equation}
  \label{eq:integration} \ell_{W,\vf}(\alpha):=\int_W\vf
  \alpha:=\int_U\vf\circ\Phi(x)\Phi^*\alpha(\Phi(x), D\Phi(x)\R^{d_s})
  \end{equation}
where $\Phi:U\to W$ is a chart and $\R^{d_s}$ is taken with the
orientation determined by corresponding elements of the
Grassmannian.\footnote{If $W$ cannot be covered by only one chart,
then the definition is trivially extended, as usual, by using a
partition of unity. Recall that, given a differential form $\omega$
on $W$ and a base $\{e_i\}$ with its dual base $\{d x_i\}$ on
$\R^{d_s}$, $\Phi^*\omega=\omega(D\Phi e_1,\dots, D\Phi
e_{d_s})\,dx_1\wedge\dots\wedge d x_{d_s}$.} A direct computation
shows that this definition is independent of the chart $\Phi$, hence
intrinsic.
\begin{rem}
\label{rem:FormsInB} If $\omega$ is a $d_s$-differential form, then
for each $(x,E)\in\G$ we can define $\alpha(x,E)$ to be the
restriction of $\omega(x)$ to $E$. Thus the forms can be embedded in
$\Se$.
\end{rem}
\begin{rem}
\label{rem:MetricsinB} A Riemannian metric defines a volume form on
any subspace of the tangent bundle of $X$. Thus, it defines an
element of $\Se$ with the property that its integral along any
nonempty compact $d_s$ dimensional submanifold is positive.
\end{rem}

\subsubsection{Integration of elements of $\Se$}
If $f:\G \to \C$ is $\Co^{r-1}$, then it is possible to multiply an
element of $\Se$ by $f$, to obtain a new element of $\Se$. In
particular, if $\alpha\in \Se$, $W$ is a $d_s$ dimensional oriented
$\Co^1$ manifold and $\vf\in\Co^0(W,\C)$, then there is a well
defined integral
  \begin{equation}
  \label{DefIntegral1}
  \int_W \vf \cdot (f \alpha).
  \end{equation}
For $x\in W$, $\tilde f(x):=f(x,T_x W)$ is a continuous function on
$W$ and so is the function $\vf \tilde f$. Hence, the integral
  \begin{equation}
  \label{DefIntegral2}
  \int_W (\vf \tilde f)\cdot  \alpha
  \end{equation}
is also well defined. By construction, the integrals
\eqref{DefIntegral1} and \eqref{DefIntegral2} coincide.
\begin{conv}
\label{Confound} We will write $\int_W \vf f \alpha$ indifferently
for these two integrals. More generally, implicitly, when we are
working along a submanifold $W$, we will confuse $f$ and $\tilde f$.
\end{conv}

\subsubsection{Lie derivative of elements of $\Se$}

If $\phi$ is a local diffeomorphism of $X$, it can be lifted through
its differential to a local bundle isomorphism of $\E$. Hence, if
$\alpha\in \Se$, its pullback $\phi^* \alpha$ is well defined. In a
pedestrian way, an element of $\Se$ is a function from $\F:=\{ (x,E,
e_1 \wedge \dots \wedge e_{d_s}) \st (x,E)\in \G, e_1,\dots, e_{d_s}
\in E\}$ to $\C$, satisfying the homogeneity relation $\alpha( x, E,
\lambda e_1\wedge \dots \wedge e_{d_s})= \lambda
\alpha(x,E,e_1\wedge\dots \wedge e_{d_s})$. If $(x,E)\in \G$ and
$e_1,\dots,e_{d_s}$ is a family of vectors in $E$, then $\phi^*
\alpha$ is given by
  \begin{equation}
  \label{eq:DefAction}
  (\phi^* \alpha)(x,E, e_1\wedge\dots \wedge e_{d_s})
  = \alpha( \phi(x), D\phi(x)E, D\phi(x)e_1 \wedge \dots
  \wedge D\phi(x) e_{d_s}).
  \end{equation}
Given a vector field $v$, we will write $L_v$ for its Lie
derivative. Given a $\Co^k$ vector field $v$ on $X$, with $k\geq 1$,
there is a canonical way to lift it to a $\Co^{k-1}$ vector field on
$\F$, as follows. Let $\phi_t$ be the flow of the vector field $v$.
For $\alpha\in \Se$, the pullback $\phi_t^* \alpha$ is well defined.
The quantity $\left.\frac{\dd \phi_t^* \alpha}{\dd t}\right|_{t=0}$
is then given by the Lie derivative of $\alpha$ against a
$\Co^{k-1}$ vector field, which we denote by $v^\F$.
 The following result will be helpful in
the following:
\begin{prop}[\cite{GeoDiffNaturelle}, Lemma 6.19]
\label{LieCommutes} The map $v\mapsto v^\F$ is linear. Moreover, if
$v_1, v_2$ are two $\Co^2$ vector fields on $X$,
  \begin{equation}
  [L_{v^\F_1},L_{v^\F_2}]=L_{[v_1,v_2]^\F}.
  \end{equation}
\end{prop}

\begin{rem}
We will use systematically the above proposition to confuse $v$ and
$v^\F$, so in the following we will suppress the superscript $\F$,
where this does not create confusion.
\end{rem}

If $W$ is a compact submanifold of $X$ with boundary, and $q\in
\R_+$, we will write $\Co^q_0(W)$ for the set of $\Co^q$ functions
from $W$ to $\C$ vanishing on the boundary of $W$, and
$\vectfield^q(W)$ for the set of $\Co^q$ vector fields defined on a
neighborhood of $W$ in $X$.

If $v\in\vectfield^1(W)$ is tangent to $W$ along $W$, and
$\alpha\in\Se$, then $L_v \alpha$ can also be obtained along $W$ by
considering the restriction of $\alpha$ to $W$, which is a volume
form, and then taking its (usual) Lie derivative with respect to the
restriction of $v$ to $W$. Therefore, the usual Stokes formula still
applies in this context, and gives the following integration by
parts formula.
\begin{prop}
Let $W$ be a compact submanifold with boundary of dimension $d_s$,
let $\alpha\in \Se$, let $v\in\vectfield^1(W)$ be tangent to $W$
along $W$, and let $\vf \in \Co^1_0(W)$. Then
  \begin{equation}
  \label{Stokes}
  \int_W \vf L_v \alpha = - \int_W (L_v \vf) \alpha.
  \end{equation}
\end{prop}

\begin{dft}
There is no divergence term. In fact, even the notion of divergence
is a ``riemannian'' notion and is not defined here. In general, if
$v$ is a vector field on some manifold, $f$ is a function and
$\alpha$ is a volume form, then $\int L_v( f\alpha)=0$ by Stokes'
Theorem (for differential forms). Since $L_v(f\alpha)=L_vf\cdot
\alpha+f\cdot L_v\alpha$, this gives
  \begin{equation}
  \int f\cdot L_v\alpha = - \int L_v f\cdot \alpha.
  \end{equation}
\end{dft}

\subsection{The norms}

Let $\Sigma$ be a set of $d_s$ dimensional compact $\Co^r$
submanifolds of $X$, with boundary. To such a $\Sigma$, we will
associate a family of norms on $\Se$ as follows.

\begin{defin}
A triple $(t,q,\coe) \in \N \times \R_+ \times \{0,1\}$ is
\emph{correct} if $t+q\leq r-1+\coe$, and $q\geq \coe$ or $t=0$.
\end{defin}

\begin{rem}
Notice that, if $(t,q,\coe)$ is correct and $t\geq 1$, then
$(t-1,q+1,\coe)$ is also correct.
\end{rem}

For any correct $(t,q,\coe)$, consider the set
  \begin{multline*}
  \Omega_{t,q+t,\coe}=\{ (W,\vf, v_1,\dots, v_t) \st W\in \Sigma,
  \vf \in \Co^{q+t}_0(W) \text{ with }|\vf|_{\Co^{q+t}(W)}\leq 1,
  \\
  v_1,\dots,v_t\in \vectfield^{q+t-\coe}(W)
  \text{ with }|v_i|_{\vectfield^{q+t-\coe}(W)}\leq 1\}.
  \end{multline*}
To each $\omega\in \Omega_{t,q+t,\coe}$, we can associate a linear
form $\ell_\omega$ on $\Se$, by
  \begin{equation}
  \ell_\omega(\alpha)= \int_W \vf \cdot L_{v_1}\dots L_{v_t}(\alpha).
  \end{equation}
Indeed, this is clearly defined if $t=0$. Moreover, if $t>0$, the
vector field $v_t$ is $\Co^{q+t-\coe}$, and is in particular
$\Co^1$. Hence, the lifted vector field $v_t^\F$ is well defined,
and $L_{v_t}(\alpha) \in \Co^{\min(r-2,
q+t-\coe-1)}=\Co^{q+t-\coe-1}$ since $q+t-\coe\leq r-1$. Going down
by induction, we have in the end $L_{v_1}\dots L_{v_t}(\alpha)\in
\Co^{q-\coe}$, which does not create any smoothness problem since
$q\geq \coe$.

We can then define the seminorms
  \begin{equation}
  \|\alpha\|_{t,q+t,\coe}^-:=\sup_{\omega\in \Omega_{t,q+t,\coe}}
  \ell_\omega(\alpha).
  \end{equation}
For $p\in \N$, $q\geq 0$ and $\coe\in \{0,1\}$ such that
$(p,q,\coe)$ is correct, we define then
  \begin{equation}
  \label{eq:pq-norm}
  \|\alpha\|_{p,q,\coe}:=\sum_{t=0}^p \|\alpha\|_{t,q+t,\coe}^-.
  \end{equation}
We will use the notation $\B^{p,q,\coe}$ for the closure of $\Se$ in
the above seminorm. This construction is as described in Section
\ref{subsec:space}.

Note that \eqref{eq:pq-norm} defines in general only a seminorm on
$\Se$. Indeed, if $\alpha \in \Se$ vanishes in a neighborhood of the
tangent spaces to elements of $\Sigma$, then
$\|\alpha\|_{p,q,\coe}=0$.

\section{The dynamics}\label{sec:due}

In Section \ref{sec:one} the dynamics did not play any role, yet all
the construction depends on the choice of $\Sigma$. In fact, such a
choice encodes in the geometry of the space the relevant properties
of the dynamics. In this chapter we will first define $\Sigma$ by
stating the relevant properties it must enjoy, then define the
transfer operator and study its properties when acting on the
resulting spaces.

\subsection{Admissible leaves}

Recall from the introduction that we have an open set $U\subset X$
and a map $T\in\Co^r(U,X)$, diffeomorphic on its image. Furthermore
$\Lambda:=\bigcap_{n\in\Z}T^n U$ is non empty and compact and
$\Lambda$ is a hyperbolic set for $T$. In addition, once and for
all, we fix an open neighborhood $U'$ of $\Lambda$, with compact
closure in $U$, such that $TU'\subset U$ and $T^{-1}U'\subset U$,
and small enough so that the restriction of $T$ to $U'$ is still
hyperbolic. For $x\in U'$, denote by $C_s(x)$ the stable cone at
$x$. Let finally $V$ be a small neighborhood of $\Lambda$, compactly
contained in $U'$.

%We will use the norms of the previous paragraph with respect to some
%set of leaves $\Sigma$. This set will have to satisfy several
%technical properties. They are all easy consequences of the
%construction of $\Sigma$ but, for definiteness and easy reference in
%the forthcoming proofs, we rather give all the properties we will need
%in the following definition.

\begin{defin}
\label{DefAdmissibleLeaves} A set $\Sigma$ of $d_s$ dimensional
compact submanifolds of $U'$ with boundary is an \emph{admissible
set of leaves} if
\begin{enumerate}
\item Each element $W$ of $\Sigma$ is a $\Co^r$ submanifold of $X$,
its tangent space at $x\in W$ is contained in $C_s(x)$, and
$\sup_{W\in \Sigma} |W|_{\Co^r} < \infty$. Moreover, for any point
$x$ of $\Lambda$, there exists $W\in \Sigma$ containing $x$ and
contained in $W^s(x)$. Additionally, $\sup_{W\in\Sigma}
\diam(W)<\infty$, and there exists $\ve>0$ such that each element of
$\Sigma$ contains a ball of radius $\ve$. Moreover, to each leaf
$W\in \Sigma$ intersecting $V$, we associate an \emph{enlargement}
$W^e$ of $W$, which is the union of a uniformly bounded number of
leaves $W_1,\dots, W_k \in \Sigma$, containing $W$, and such that
$\dist(\partial W,\partial W^e)>2\delta_0$ for some $\delta_0>0$
(independent of $W$).

\item Let us say that two leaves $W, W'\in \Sigma$ are $(C,\ve)$-close if
there exists a $\Co^{r-1}$ vector field $v$, defined on a
neighborhood of $W$, with $|v|_{\Co^{r-1}} \leq \ve$, and such that
its flow $\phi_t$ is uniformly bounded in $\Co^r$ by $C$ and
satisfies $\phi_1(W)=W'$ and $\phi_t(W)\in\Sigma$ for $0\leq t\leq
1$. We assume that there exists a constant $C_\Sigma$ such that, for
all $\ve>0$, there exists a finite number of leaves $W_1,\dots,
W_k\in \Sigma$ such that any $W\in \Sigma$ is $(C_\Sigma,\ve)$-close
to a leaf $W_i$ with $1\leq i\leq k$.

\item There exist $C>0$ and a sequence $\ve_n$ going
exponentially fast to $0$ such that, for all $n\in \N^*$, for all
$W\in \Sigma$, there exist a finite number of leaves
$W_1,\dots,W_k\in \Sigma$ and $\Co^r$ functions $\rho_1,\dots,
\rho_k$ with values in $[0,1]$ compactly supported on $W_i$, with
$|\rho_i|_{\Co^r(W_i)} \leq C$, and such that the set
$W^{(n)}=\{x\in W, \forall\; 0\leq i\leq n-1, T^{-i}x \in V\}$
satisfies: $T^{-n}W^{(n)} \subset \bigcup W_i$, and $\sum \rho_i=1$
on $T^{-n}W^{(n)}$, and any point of $T^{-n}W^{(n)}$ is contained in
at most $C$ sets $W_i$. Moreover, $W_i$ is $(C_\Sigma,\ve_n)$-close
to an element of $\Sigma$ contained in the stable manifold of a
point of $\Lambda$. Finally, $T^n\left(\bigcup_{i=1}^k W_i\right)$
is contained in the enlargement $W^e$ of $W$, and even in the set
$\{ x\in W^e \st \dist(x,\partial W^e)>\delta_0\}$.
\end{enumerate}
\end{defin}
The first property of the definition means that the elements of
$\Sigma$ are close to stable leaves in the $\Co^1$ topology, and
have a reasonable size.
 The second condition means that there are
sufficiently many leaves, and will imply some compactness
properties. The third property is an invariance property and means
that we can iterate the leaves backward.

In \cite{gouezel_liverani}, the existence of admissible sets of
leaves is proved for Anosov systems. The proof generalizes in a
straightforward way to this setting. Hence, the following
proposition holds.
\begin{prop}
Admissible sets of leaves do exist.
\end{prop}

We choose once and for all such an admissible set of leaves, and
denote it by $\Sigma$.

\subsection{Definition of the Operator}

\label{subsec:transfer} We will consider the action of the
composition by $T$ on the previously defined spaces. For historical
reasons, we will rather consider the composition by $T^{-1}$, but
this choice is arbitrary. To keep the functions supported in $U$, we
need a truncation function. Let $\trunc$ be a $\Co^r$ function
taking values in $[0,1]$, equal to $1$ on a neighborhood of
$\Lambda$ and compactly supported in $T(V)$.

We need also to introduce a \emph{weight}. We will consider two
classes of weights. Let $\Weight^0$ be the set of $\Co^{r-1}$
functions $\phi$ from $\G$ to $\R$, such that if $x\in U$ and $F$
and $F'$ are the same subspace of $T_x U$ but with opposite
orientations, then $\phi(F)=\phi(F')$. This condition makes it
possible to define a function $\bar \phi$ on $\Lambda$ by $\bar
\phi(x)=\phi(x, E^s(x))$ (where the orientation of $E^s(x)$ is not
relevant by the previous property). Let $\Weight^1$ be the set of
$\Co^r$ functions from $X$ to $\R$. Of course, an element of
$\Weight^1$ is an element of $\Weight^0$ as well. Yet, slightly
stronger results hold true for weights in $\Weight^1$.

For each truncation function $\trunc$, and each weight
$\phi\in\Weight^\coe$, $\coe\in \{0,1\}$, we define a
\emph{truncated and weighted transfer operator} (or simply transfer
operator) $\Lp_{\trunc, \phi}:\Se\to\Se$ by
  \begin{equation}
  \Lp_{\trunc, \phi}\alpha(x,E):=\trunc(T^{-1}x) e^{ \phi(T^{-1}x,
  DT^{-1}(x)E)} T_*\alpha (T^{-1}x, DT^{-1}(x)E).
  \end{equation}
In terms of the action of diffeomorphisms on elements of $\Se$
defined in \eqref{eq:DefAction}, this formula can be written as
$\Lp_{\trunc,\phi} \alpha=T_*(\pi e^\phi \alpha)$. It is clear that
an understanding of the iterates of $\Lp_{\trunc,\phi}$ would shed
light on the mixing properties of $T$. The operator
$\Lp_{\trunc,\phi}$ does not have good asymptotic properties on
$\Se$ with its $\Co^{r-1}$ norm, but we will show that it behaves
well on the spaces $\B^{p,q,\coe}$.

If $W$ is a submanifold of dimension $d_s$ contained in $U'$, $\vf$
is a continuous function on $W$ with compact support and $\alpha\in
\Se$, then by definition
  \begin{equation}\label{eq:transfer}
  \int_W \vf \Lp_{\trunc, \phi}\alpha= \int_{T^{-1}W} \vf\circ T
  \trunc
  e^{\phi} \alpha.
  \end{equation}
Recall that this integral is well defined by Convention
\ref{Confound}.

\subsection{Main dynamical inequality}

When $(p,q,\coe)$ is correct, i.e., $p+q-\coe\leq r-1$, and $q\geq
\coe$ or $t=0$, and the weight $\phi$ belongs to $\Weight^\coe$, we
can study the spectral properties of $\Lp_{\trunc, \phi}$ acting on
$\B^{p,q,\coe}$. Notice that, for a weight belonging to $\Weight^1$,
this means that we can go up to $p+q=r$, i.e., we can reach the
differentiability of the map. Before proceeding we need a
definition.
\begin{defin}
\label{def:radius} For each $W\in\Sigma$ and $n\in\N$ let $\{W_j\}$
be any covering of $T^{-n}W^{(n)}$ as given by the third item in
Definition \ref{DefAdmissibleLeaves}. We define
  \begin{equation}
  \ra_n:=
  \left(\sup_{W\in\Sigma}\sum_j\left|e^{S_n \phi} \trunc_n
  \right|_{\Co^{r-1+\coe}(W_j)}\right)^{1/n},
  \end{equation}
where $\trunc_n:=\prod_{k=0}^{n-1}\trunc\circ T^k$ and, for each
function $f:X\to\C$, $S_n f:=\sum_{k=0}^{n-1}f\circ T^k$. Here, to
define $S_n \phi$ along $W$, we use Convention
\ref{Confound}.\footnote{\label{foot:exist}Note that the volume of
$T^{-n}W$ grows at most exponentially. Thus, given the condition (3)
of Definition \ref{DefAdmissibleLeaves} on the bounded overlap of
the $W_j$, the cardinality of $\{W_j\}$ can grow at most
exponentially as well. In turn, this means that there exists a
constant $C$ such that $\ra_n \leq C$.}
\end{defin}

The main lemma to prove Lasota-Yorke type inequalities is the
following:
\begin{lem}
\label{MainDynamicalInequality} Let $t\in \N$ and $q\geq 0$. If
$(t,q,\coe)$ is correct, there exist constants $C>0$ and $C_n>0$ for
$n\in \N$ such that, for any $\alpha\in \Se$,
  \begin{equation}
  \label{LY:eq2}
  \norm{\Lp_{\trunc,\phi}^n \alpha}_{t,q+t,\coe}^-
  \leq C \ra_n^n \lambda^{-tn}
  \norm{\alpha}_{t,q+t,\coe}^- + C_n
  \sum_{0\leq t'<t} \norm{\alpha}_{t',q+t',\coe}^-.
  \end{equation}
Moreover, if $(t,q+1,\coe)$ is also correct,
  \begin{multline}
  \label{LY:eq1}
  \norm{\Lp^n_{\trunc,\phi} \alpha}_{t,q+t,\coe}^- \leq C \ra_n^n
  \nu^{(q+t)n} \lambda^{-tn}
  \norm{\alpha}_{t,q+t,\coe}^-
  \\+C_n\sum_{0\leq t'<t} \norm{\alpha}_{t',q+t',\coe}^-
  + C_n \norm{\alpha}_{t,q+t+1,\coe}^-.
  \end{multline}
\end{lem}
\begin{proof}
Take $\omega=(W,\vf, v_1,\dots,v_t)\in \Omega_{t,q+t,\coe}$. Let
$\rho_j$ be an adapted partition of unity on $T^{-n}W^{(n)}$, as
given in (3) of Definition \ref{DefAdmissibleLeaves}. We want to
estimate
  \begin{equation}
  \label{to_estimate}
  \int_W \vf \cdot L_{v_1}\dots L_{v_t}( \Lp_{\trunc,\phi}^n \alpha)
  = \sum_j \int_{W_j} \vf\circ T^n \rho_j
  \cdot L_{w_1}\dots L_{w_t}( \alpha \cdot
  e^{S_n \phi}\prod_{k=0}^{n-1} \pi\circ
  T^k ),
  \end{equation}
where $w_i=(T^n)^*(v_i)$. Remembering that $\pi_n:=\prod_{k=0}^{n-1}
\pi\circ
  T^k$,
  \begin{equation}
  \label{decompose}
  L_{w_1}\dots L_{w_t}(\alpha \cdot e^{S_n \phi}\pi_n)= \sum_{A\subset
  \{1,\dots, t\}} \left(\prod_{i\not\in A} L_{w_i}\right) (\alpha) \cdot
  \left(\prod_{i\in A} L_{w_i}\right)(e^{S_n \phi}\pi_n).
  \end{equation}
We claim that, for any $A\subset \{1,\dots, t\}$,
  \begin{equation}
  \label{belongs_good_set}
  \left( \prod_{i
  \in A} L_{w_i}\right)(e^{S_n \phi}\pi_n) \in \Co^{q+t-\# A}.
  \end{equation}
Assume first that $\coe=0$. Then $q+t \leq r-1$. The lift $w_k^\F$
of any of the vector fields $w_k$ is in $\Co^{q+t-1}$, hence
$L_{w_k}(e^{S_n \phi}\pi_n) \in \Co^{\min(r-2, q+t-1)}=\Co^{q+t-1}$.
Equation \eqref{belongs_good_set} then follows inductively on $\#
A$. On the other hand, if $\coe=1$, the vector field $w_k$ is only
$\Co^{q+t-1}$ (and so $w_k^\F$ is only $\Co^{q+t-2}$, which is not
sufficient). However, there is no need to lift the vector field
$w_k$ to $\F$ since $\phi$ is defined on $X$. Hence, we get
$L_{w_k}(e^{S_n \phi}\pi_n) \in \Co^{\min(r-1,q+t-1)}=\Co^{q+t-1}$.
Equation \eqref{belongs_good_set} easily follows.

In the right hand side of \eqref{to_estimate}, we can use
\eqref{decompose} to compute $L_{w_1}\dots L_{w_t}( \alpha \cdot
  e^{S_n \phi}\pi_n)$. Any term with $A\not=\emptyset$ is then estimated
by $C_n \norm{\alpha}^-_{t-\# A, q+t-\#A, \coe}$, thanks to
\eqref{belongs_good_set}. Hence, to conclude, it suffices to
estimate the remaining term with $A=\emptyset$\ ~:
  \begin{equation}
  \int_{W_j} \vf\circ T^n \rho_j e^{S_n \phi}
  \pi_n \cdot L_{w_1}\dots L_{w_t} (
  \alpha).
  \end{equation}
To this end we decompose $w_i$ as $w_i^u+w_i^s$ where $w_i^u$ and
$w_i^s$ are $\Co^{q+t-\coe}$ vector fields, $w_i^s$ is tangent to
$W_j$, and $|w_i^u|_{\Co^{q+t-\coe}} \leq C
\lambda^{-n}$.\footnote{\label{ft:decompose}Such a
decomposition is achieved in  \cite[Appendix A]{gouezel_liverani}
(and the computation is even easier since the smoothing is not
required). The argument roughly goes as follows. Consider a
$\Co^{r}$ foliation transverse to $W_j$, and push it by $T^n$.
Around $T^n W_j$, consider also a foliation given by translates (in
some chart with uniformly bounded $\Co^{r}$ norm) of $T^n W_j$. Then
project simply $v_i$ on these two transverse foliations, and pull
everything back under $T^n$. This is essentially the desired
decomposition.} Clearly
$L_{w_i}=L_{w_i^u}+L_{w_i^s}$. Hence, for $\sigma\in \{s,u\}^t$, we
must study the integrals
  \begin{equation}
  \label{to_study}
  \int_{W_j} \vf\circ T^n \rho_j e^{S_n \phi} \pi_n\cdot
  L_{w^{\sigma_1}_1} \dots L_{w^{\sigma_t}_t} (
  \alpha).
  \end{equation}
Notice that, if we exchange two of these vector fields, the
difference is of the form $\int_{W_j} \tilde\vf L_{w'_1}\dots
L_{w'_{t-1}}(\alpha)$ where $w'_1,\dots, w'_{t-1}$ are
$\Co^{q+t-1-\coe}$ vector fields. Indeed, $L_w L_{w'}= L_{w'}L_w
+L_{[w,w']}$ by Proposition \ref{LieCommutes}, and $[w,w']$ is a
$\Co^{q+t-1-\coe}$ vector field. In particular, up to $C_n
\norm{\alpha}_{t-1, q+t-1,\coe}^-$, we can freely exchange the
vector fields.

Suppose first that $\sigma_1=s$. Then, by \eqref{Stokes}, the
integral \eqref{to_study} is equal to
  \begin{equation}
  - \int_{W_j} L_{w_1^s}(  \vf\circ T^n \rho_j e^{S_n
  \phi}\pi_n) \cdot
  L_{w^{\sigma_2}_2} \dots L_{w^{\sigma_t}_t} (
  \alpha).
  \end{equation}
This is bounded by $C_n \norm{\alpha}_{t-1, q+t-1,\coe}^-$. More
generally, if one of the $\sigma_i$'s is equal to $s$, we can first
exchange the vector fields as described above to put the
corresponding $L_{w_i^s}$ in the first place, and then integrate by
parts. Finally, we have
  \begin{multline}
  \label{eq:jla;ksjfd}
  \int_{W_j} \vf\circ T^n \rho_j e^{S_n \phi} \pi_n\cdot L_{w_1}\dots L_{w_t} (
  \alpha)\\=
  \int_{W_j} \vf\circ T^n \rho_j e^{S_n \phi} \pi_n\cdot L_{w^u_1}\dots L_{w^u_t} (
  \alpha) + \Or( \norm{\alpha}_{t-1,q+t-1,\coe}^-).
  \end{multline}

We are now positioned to prove \eqref{LY:eq2}. The last integral in
\eqref{eq:jla;ksjfd} is bounded by
  \begin{multline}
  \label{eq:tequalp}
  \left| \vf\circ T^n \rho_j e^{S_n
  \phi} \pi_n \right|_{\Co^{q+t}(W_j)}
  \prod_{i=1}^t |w^u_i|_{\Co^{q+t-\coe}(W_j)}
  \norm{\alpha}_{t,q+t,\coe}^-
  \\
  \leq C
  \lambda^{-tn} \left| e^{S_n \phi}\pi_n\right|_{\Co^{q+t}(W_j)}
  \norm{\alpha}_{t,q+t,\coe}^-.
  \end{multline}
Summing the inequalities \eqref{eq:tequalp} over $j$ and remembering
Definition \ref{def:radius} yields \eqref{LY:eq2}.

This simple argument is not sufficient to prove \eqref{LY:eq1},
since we want also to gain a factor $\nu^{-(q+t)n}$ (if we are ready
to pay the price of having a term $\norm{\alpha}_{t, q+t+1,\coe}^-$
in the upper bound). To do this, we will smoothen the test function
$\vf$. Let $\A_\ve \vf$ be obtained by convolving $\vf$ with a
mollifier of size $\ve$. If $a$ is the largest integer less than
$q+t$, we have $|\A_\ve \vf - \vf|_{\Co^a} \leq C \ve^{q+t-a}$, the
function $\A_\ve \vf$ is bounded in $\Co^{q+t}$ independently of
$\ve$, and it belongs to $\Co^{q+t+1}$. We choose $\ve=\nu^{(q+t)n/
(q+t-a)}$. In this way,
  \begin{equation}
  | (\vf - \A_\ve \vf)\circ T^n |_{\Co^{q+t}(W_j)} \leq C
  \nu^{(q+t)n}.
  \end{equation}
Then \eqref{eq:jla;ksjfd} implies
  \begin{align*}
  \int_{W_j} \vf\circ T^n \rho_j e^{S_n \phi} \pi_n
  \cdot L_{w^u_1}\dots L_{w^u_t} (
  \alpha)
  \!\!\!\!\!\!\!\!\! \!\!\!\!\!\!\!\!\! \!\!\!\!\!\!\!\!\!
   \!\!\!\!\!\!\!\!\! &
  \\&
  = \int_{W_j} (\vf- \A_\ve \vf) \circ T^n \rho_j e^{S_n \phi}
  \pi_n
  \cdot L_{w^u_1}\dots L_{w^u_t} ( \alpha)
  \\&
  + \int_{W_j} (\A_\ve \vf)\circ T^n
  \rho_j e^{S_n \phi} \pi_n
  \cdot L_{w^u_1}\dots L_{w^u_t} (\alpha).
  \end{align*}
The last integral is bounded by $C_n \norm{\alpha}_{t,
q+t+1,\coe}^-$. And the previous one is at most
  \begin{multline*}
  \left|  (\vf- \A_\ve \vf) \circ T^n \rho_j e^{S_n
  \phi}\pi_n\right|_{\Co^{q+t}(W_j)} \prod_{i=1}^t |
  w^u_i|_{\Co^{q+t-\coe}(W_j)}\norm{\alpha}_{t,q+t,\coe}^-
  \\
  \leq C \nu^{(q+t)n} \left|e^{S_n \phi}\pi_n\right|_{\Co^{q+t}(W_j)}
  \lambda^{-tn} \norm{\alpha}_{t,q+t,\coe}^-.
  \end{multline*}
Summing over $j$ and remembering Definition \ref{def:radius}, we
finally have \eqref{LY:eq1}.
\end{proof}

\section{Spectral properties of the Transfer Operator}
\label{sec:specral-radius} In this section we investigate the
spectral radius and the essential spectral radius of the Ruelle
operator. We will use constants $\ora>0$ and $d\in \N$ such
that\footnote{Such constants do exist, see footnote
\ref{foot:exist}.}
  \begin{equation}
  \label{eq:defora}
  \exists C>0,\forall n\in \N^*,\quad \ra_n^n \leq C n^d \ora^n.
  \end{equation}

\subsection{Quasi compactness}

As usual, the proof of the quasi compactness of the transfer
operator is based on two ingredients: a Lasota-Yorke type inequality
and a compact embedding between spaces. See \cite{baladi:decay} if
unfamiliar with such ideas.

\subsubsection{Lasota-Yorke inequality}

\begin{lem}
\label{lem:LY} Let $\coe\in \{0,1\}$. For all $p\in \N$ and $q\geq
0$ such that $(p,q,\coe)$ is correct, for all $(\ora,d)$ satisfying
\eqref{eq:defora}, there exists a constant $C>0$ such that, for all
$n\in \N^*$, for all $\alpha\in \Se$,
  \begin{equation}
  \label{eq:bounded}
  \norm{\Lp_{\trunc,\phi}^n \alpha}_{p,q,\coe} \leq C n^d \ora^n
  \norm{\alpha}_{p,q,\coe}.
  \end{equation}
Moreover, if $p>0$, the following inequality also holds:
  \begin{equation}
  \label{eq:bounded2}
  \norm{\Lp_{\trunc,\phi}^n \alpha}_{p,q,\coe} \leq C \ora^n
  \max(\lambda^{-p},\nu^q)^n \norm{\alpha}_{p,q,\coe} + C n^d \ora^n
  \norm{\alpha}_{p-1,q+1,\coe}.
  \end{equation}
Finally, if $p\geq 0$, there exists $\sigma<1$ (independent of
$\ora$ and $d$) such that
  \begin{equation}
  \label{eq:bounded3}
  \norm{\Lp_{\trunc,\phi}^n \alpha}_{p,q,\coe} \leq C \ora^n
  \sigma^n \norm{\alpha}_{p,q,\coe} + C n^d \ora^n
  \norm{\alpha}_{0,q+p,\coe}.
  \end{equation}
\end{lem}
\begin{proof}
The inequality \eqref{eq:bounded3} is an easy consequence of
\eqref{eq:bounded2} and an induction on $p$. Moreover,
\eqref{eq:bounded} for $p=0$ is a direct consequence of Equation
\eqref{LY:eq2} with $t=0$, and \eqref{eq:defora}. Note also that
\eqref{eq:bounded2} for $p>0$ implies \eqref{eq:bounded} for the
same $p$. Hence, it is sufficient to prove that \eqref{eq:bounded}
at $p-1$ implies \eqref{eq:bounded2} at $p$.

Choose any $\lambda'>\lambda$ and $\nu'<\nu$ respectively smaller
and larger than the best expansion and contraction constants of $T$
in the unstable and stable direction. Lemma
\ref{MainDynamicalInequality} still applies with $\lambda'$ and
$\mu'$ instead of $\lambda$ and $\mu$. Hence, there exist constants
$C_0$ and $C'_n$ such that, for all $0\leq t \leq p$, and setting
$\sigma_1:=\max({\lambda'}^{-p}, {\nu'}^q)$,
\[
  \norm{\Lp_{\trunc,\phi}^n \alpha}_{t,q+t,\coe}^-
  \leq C_0 n^d \ora^n \sigma_1^n
  \norm{\alpha}_{t,q+t,\coe}^-
  + C'_n \sum_{t'<t}
  \norm{\alpha}_{t',q+t',\coe}^- + C'_n \norm{\alpha}_{p-1,q+1,\coe}.
\]
To prove this, we use \eqref{LY:eq2} for $t=p$, and \eqref{LY:eq1}
for $t<p$ (in which case $\norm{\alpha}_{t,q+t+1,\coe}^- \leq
\norm{\alpha}_{p-1,q+1,\coe}$).

Let $\sigma= \max(\lambda^{-p}, \nu^q)$. There exists $N$ such that
$C_0 N^d \sigma_1^N \leq \sigma^N /2$. We fix it once and for all.
Fix also once and for all a large constant $K>2$ such that
$\frac{C'_N K^{-1}}{1-K^{-1}} \leq \ora^N \sigma^N/2$, and define a
new seminorm on $\Se$ by $\norm{\alpha}_{p,q,\coe}'=\sum_{t=0}^p
K^{-t} \norm{\alpha}_{t,q+t,\coe}^-$. Then
$\norm{\Lp_{\trunc,\phi}^N \alpha}_{p,q,\coe}'$ is at most
  \begin{align*}
  &\sum_{t=0}^p K^{-t}
  \Biggl( \ora^N (\sigma^N/2)
  \norm{\alpha}_{t,q+t,\coe}^- +
  C'_N \sum_{t'<t} \norm{\alpha}_{t',q+t',\coe}^-
  + C'_N \norm{\alpha}_{p-1,q+1,\coe} \Biggr)
  \\&
  \leq \ora^N(\sigma^N /2) \norm{\alpha}_{p,q,\coe}' + C'_N \sum_{t'=0}^p
  \frac{K^{-t'-1}}{1-K^{-1}} \norm{\alpha}_{t', q+t',\coe}^- + \frac{C'_N}{1-K^{-1}}
  \norm{\alpha}_{p-1, q+1,\coe}
  \\&
  \leq \ora^N
  (\sigma^N /2) \norm{\alpha}_{p,q,\coe}' + \frac{C'_N K^{-1}}{1-K^{-1}}
  \norm{\alpha}_{p,q,\coe}' +2C'_N \norm{\alpha}_{p-1, q+1,\coe}.
  \end{align*}
Since $K$ was chosen large enough, we have therefore
  \begin{equation}
  \norm{\Lp_{\trunc,\phi}^N \alpha}_{p,q,\coe}' \leq
  \ora^N \sigma^N \norm{\alpha}_{p,q,\coe}'+ 2C'_N
  \norm{\alpha}_{p-1, q+1,\coe}.
  \end{equation}
By the inductive assumption, the iterates of $\Lp_{\trunc,\phi}$
satisfy the inequality
  \begin{equation}
  \norm{\Lp_{\trunc,\phi}^n \alpha}_{p-1,q+1,\coe}\leq C_1 n^d \ora^n
  \norm{\alpha}_{p-1,q+1,\coe},
  \end{equation}
for some constant $C_1$. This implies by induction on $m$ that
  \begin{align*}
  \norm{\Lp_{\trunc,\phi}^{mN} \alpha}_{p,q,\coe}'
  &
  \leq (\ora\sigma)^{mN} \norm{\alpha}_{p,q,\coe}'
  + 2C'_N \sum_{k=1}^m (\ora \sigma)^{(k-1)N} \norm{
  \Lp_{\trunc,\phi}^{(m-k)N}\alpha}_{p-1,q+1,\coe}
  %
  %\\&
  %\leq (\ora\sigma)^{mN} \norm{\alpha}_{p,q,\coe}'
  %+2C'_N\sum_{k=1}^m (\ora \sigma)^{(k-1)N}C_1 (mN)^d \ora^{(m-k)N}
  %\norm{\alpha}_{p-1,q+1,\coe}
  %
  \\&
  \leq (\ora\sigma)^{mN} \norm{\alpha}_{p,q,\coe}'
  +2C'_N C_1 (mN)^d \ora^{mN} \ora^{-N}
  \left( \sum_{i=0}^\infty \sigma^{iN} \right)
  \norm{\alpha}_{p-1,q+1,\coe}.
  \end{align*}
Finally, taking care of the first $N$ iterates, we obtain:
  \begin{equation}
  \norm{\Lp_{\trunc,\phi}^n \alpha}_{p,q,\coe}'
  \leq C \ora^n\sigma^n \norm{\alpha}_{p,q,\coe}' + C n^d \ora^n
  \norm{\alpha}_{p-1, q+1,\coe}.
  \end{equation}
Since the norms $\norm{\cdot}_{p,q,\coe}$ and
$\norm{\cdot}_{p,q,\coe}'$ are equivalent, this concludes the proof.
\end{proof}

\subsubsection{Compact embedding of $\B^{p,q,\coe}$ in $\B^{p-1,q+1,\coe}$}

\begin{lem}
\label{lem:CloseLeaves} Assume that $(t,q,\coe)$ is correct and that
$(t+1,q-1,\coe)$ is also correct. There exists a constant $C>0$ such
that, for all $\ve>0$, for all $W,W'$ which are
$(C_\Sigma,\ve)$-close,\footnote{See Definition
\ref{DefAdmissibleLeaves} for the definition of
$(C_\Sigma,\ve)$-close.} for all $\alpha\in \Se$,
  \begin{equation*}
  \sup_{ \omega'=(W',\vf',v'_1,\dots,v'_t)\in \Omega_{t,q+t,\coe}}
  |\ell_{\omega'}(\alpha)|
  \leq C \sup_{ \omega=(W,\vf,v_1,\dots,v_t)\in \Omega_{t,q+t,\coe}}
  |\ell_{\omega}(\alpha)|+C \ve\norm{\alpha}^-_{t+1,q+t,\coe}.
  \end{equation*}
\end{lem}
\begin{proof}
Let $v$ be a vector field with $|v|_{\Co^{r-1}}\leq \ve$ whose flow
$\phi_u$ satisfies $\phi_1(W)=W'$ and $W^u=\phi_u(W)\in \Sigma$ for
$0\leq u\leq 1$, and is bounded in $\Co^r$ by $C_\Sigma$. Start from
$\omega'=(W',\vf',v'_1,\dots,v'_t)\in \Omega_{t,q+t,\coe}$. Define
vector fields $v_i^u= \phi_{1-u}^* v'_i$, $v^u=\phi_{1-u}^* v$ and
functions $\vf^u= \vf'\circ \phi_{1-u}$. Let
  \begin{equation}
  F(u)=\int_W \vf^0 L_{v_1^0}\dots L_{v_t^0} ( \phi_u^* \alpha)
  = \int_{W^u} \vf^u L_{v_1^u}\dots L_{v_t^u} ( \alpha).
  \end{equation}
Then $F(1)= \ell_{\omega'}(\alpha)$, and $F(0)=\int_W \vf^0
L_{v_1^0}\dots L_{v_t^0} (\alpha)$. Since the vector fields $v_i^0$
have a uniformly bounded $\Co^{q+t-\coe}$ norm, and $\vf^0$ is
uniformly bounded in $\Co^{q+t}$, it is sufficient to prove that
$|F(1)-F(0)| \leq C \ve\norm{\alpha}^-_{t+1,q+t,\coe}$ to conclude.
We will prove such an estimate for $F'(u)$.

We have
  \begin{equation}
  F'(u)= \int_W \vf^0 L_{v_1^0}\dots L_{v_t^0} ( \phi_u^* L_v \alpha)
  = \int_{W^u} \vf^u L_{v_1^u} \dots L_{v_t^u} L_{v^u} \alpha.
  \end{equation}
By definition of $\norm{\cdot}^-_{t+1,q+t,\coe}$, this quantity is
bounded by $C\norm{\alpha}^-_{t+1,q+t,\coe}$, which concludes the
proof.
\end{proof}

Assume that $(p,q,\coe)$ is correct and $p>0$. Hence,
$(p-1,q+1,\coe)$ is also correct. Moreover, for any $\alpha\in \Se$,
$\norm{\alpha}_{p-1,q+1,\coe} \leq \norm{\alpha}_{p,q,\coe}$. Hence,
there exists a canonical map $\B^{p,q,\coe} \to \B^{p-1,q+1,\coe}$
extending the identity on the dense subset $\Se$ of $\B^{p,q,\coe}$.
\begin{lem}
\label{lem:compactness} If $(p,q,\coe)$ is correct and $p>0$, the
canonical map from $\B^{p,q,\coe}$ to $\B^{p-1,q+1,\coe}$ is
compact.
\end{lem}
\begin{proof}
The main point of the proof of Lemma \ref{lem:compactness} is to be
able to work only with a finite number of leaves. This is ensured by
Lemma \ref{lem:CloseLeaves}. The rest of the proof is then very
similar to \cite[Proof of Lemma 2.1]{gouezel_liverani}.
\end{proof}

\subsubsection{Spectral gap}

Lemmas \ref{lem:LY} and \ref{lem:compactness}, giving a Lasota-Yorke
inequality and compactness, imply a  precise spectral description of
the transfer operator $\Lp_{\trunc,\phi}$. Let
  \begin{equation}
  \label{eq:spectral-radius}
  \ra:=\limsup_{n\to\infty}\ra_n.
  \end{equation}
\begin{prop}
\label{thm:SpectralGap} Assume that $(p,q,\coe)$ is correct. The
operator $\Lp_{\trunc,\phi}: \Se \to \Se$ extends to a continuous
operator on $\B^{p,q,\coe}$. Its spectral radius is at most $\ra$
and its essential spectral radius is at most
$\max(\lambda^{-p},\nu^q) \ra$.
\end{prop}
\begin{proof}
For any $\ora>\ra$, the inequality \eqref{eq:bounded2}, the
compactness Lemma \ref{lem:compactness} and Hennion's Theorem
\cite{hennion} prove that the spectral radius of $\Lp_{\trunc,\phi}$
acting on $\B^{p,q,\coe}$ is bounded by $\ora$, and that its
essential spectral radius is bounded by $\max(\lambda^{-p},\nu^q)
\ora$. Letting $\ora$ tend to $\ra$, we obtain the required upper
bounds on the spectral radius and essential spectral radius of
$\Lp_{\trunc,\phi}$.
\end{proof}

\subsection{A lower bound for the spectral radius}
\label{sub:lowerbound}

We will prove that the spectral radius of $\Lp_{\trunc,\phi}$ is in
fact \emph{equal} to $\ra$. To do this, we will need the following
lower bound on $\ra_n$. Since we will use this lemma again later, to
exclude the possibility of Jordan blocks, we formulate it in greater
generality than currently needed.

\begin{lem}
\label{lem:Boundsra} Assume that $(p,q,\coe)$ is correct. Let
$\alpha$ be an element of $\B^{p,q,\coe}$ which induces a
nonnegative measure on every admissible leaf $W\in \Sigma$. Assume
moreover that there exists an open set $O$ containing $\Lambda$ such
that, for any $\ve>0$, there exists $c_\ve>0$ such that, for any
$x\in O\cap \bigcap_{n\geq 0} T^n V$, for any $W\in \Sigma$
containing $x$ with $\dist(x,\partial W)>\ve$, holds
$\int_{B_W(x,\ve)} \alpha \geq c_\ve$.\footnote{Here, $B_W(x,\ve)$
denotes the ball of center $x$ and radius $\ve$ in the manifold
$W$.} Then there exist $L\in \N$ and $C>0$ such that, for all large
enough $n$,
  \begin{equation}
  \ra_n^n \leq C \norm{ \Lp_{\trunc,\phi}^{n-2L} \alpha}_{p,q,\coe}.
  \end{equation}
\end{lem}
\begin{proof}
Let $W \in \Sigma$, and let $W_j$ be a covering of $T^{-n}W^{(n)}$
as given by Definition \ref{DefAdmissibleLeaves}. All is needed is
to prove the inequality
  \begin{equation}
  \sum_j |e^{S_n \phi} \pi_n|_{\Co^{r-1+\coe}(W_j)} \leq
  C \norm{ \Lp_{\trunc,\phi}^{n-2L} \alpha}_{p,q,\coe}.
  \end{equation}
The lemma would have a two lines proof if we could use distortion to
estimate $|e^{S_n \phi} \pi_n|_{\Co^{r-1+\coe}(W_j)}$ by $\int_{W_j}
e^{S_n \phi} \pi_n \alpha$, but there are two problems in doing so.
First, $\pi$ vanishes at some points, hence classical distortion
controls do not apply. Second, the behavior of $\alpha$ is known
only for leaves close to $\Lambda$. To overcome these two problems,
we will consider a small neighborhood of $\Lambda$, where $\pi_n$ is
equal to $1$ and $\alpha$ is well behaved. We can assume without
loss of generality that $\trunc =1 $ on $O$.

Recall the definition of the constant $\delta_0$ in the first item
of Definition \ref{DefAdmissibleLeaves}. Decreasing $\delta_0$ if
necessary, we can assume that, for all $x\in \Lambda$,
$B(x,3\delta_0) \subset O$. Then there exist $\ve>0$ and a small
neighborhood $O'$ of $\Lambda$ with the following property:
\emph{let $x\in O'$, and let $Z$ be a submanifold of dimension $d_s$
containing $x$, whose tangent space is everywhere contained in the
stable cone, and with $\dist(x,\partial Z) \geq \delta_0$. Then
there exists a point $y\in Z \cap O \cap \bigcap_{n\geq 0} T^n V$
such that $\dist(y,\partial Z) \geq \ve$  and $\dist(x,y)\leq
\delta_0$.} This is a consequence of the compactness of $\Lambda$
and the uniform transversality between the stable cones and the
unstable leaves. Decreasing $O'$ if necessary, we can assume that
  \begin{equation}
  \label{GoodInclusions}
  \forall x\in O', \quad
  B(x, 2\delta_0) \subset O.
  \end{equation}
We can also assume $\ve < \delta_0$.

We will use the following fact: \emph{there exists $L\in \N$ such
that, for any point $x$, for any $n\geq 2L$, if $T^i x\in V$ for all
$0\leq i \leq n-1$ then $T^i x\in O'$ for all $L\leq i \leq n-L$.}

This is a classical property of locally maximal sets, proved as
follows. If the fact were not true, we would have for all $L\geq 0$
a point $x_L \in V\backslash O'$ such that $T^i x_L \in V$ for all
$|i|\leq L$. An accumulation point of the sequence $x_L$ would then
belong to $\overline{V} \backslash O'$, and also to $\bigcap_{n\in
\Z} T^{-n}U$. This is a contradiction since this last intersection
is equal to $\Lambda$ by assumption, and is therefore contained in
$O'$.

Let us now return to the proof. We start from the covering $\{W_j\}$
of $T^{-n}W^{(n)}$. Fix some $j$ such that $\pi_n$ is not zero on
$W_j$. There exists $x_j\in W_j$ such that $T^i x_j\in V$ for $0\leq
i \leq n-1$. The above fact ensures that $T^i x_j\in O'$ for $L\leq
i \leq n-L$. By definition of the enlargement $W^e$ of $W$, the
point $T^n x_j$ belongs to $\{ y\in W^e \st \dist(y,\partial
W^e)\geq \delta_0\}$. Since $T^{-1}$ expands the distances in the
stable cone, we get $\dist(T^L x_j, \partial( T^{-(n-L)}W^e)) \geq
\delta_0$. Therefore, the above property shows the existence of a
point $y_j \in T^{-(n-L)}W^e \cap O \cap \bigcap_{n\geq 0} T^n V$,
with  $\dist(T^L x_j, y_j) \leq \delta_0$, such that the ball $B_j$
of center $y_j$ and radius $\ve$ in the manifold $T^{-(n-L)} W^e$ is
well defined. This ball satisfies $\int_{B_j} \alpha \geq c_\ve$ by
the assumption of the lemma. Moreover, by contraction of the
iterates of $T$ along $T^{-(n-L)}W^e$, we have $T^i(B_j) \subset
B(T^{L+i} x_j, 2\delta_0)$ for $0\leq i \leq n-L$. Since $T^{L+i}
x_j \in O'$ for $0\leq i \leq n-2L$, \eqref{GoodInclusions} shows
that $T^i(B_j) \subset O$ for $0\leq i \leq n-2L$. Therefore,
$\trunc_{n-2L}=1$ on $B_j$.

By uniform contraction of $T$, $|\pi_n|_{\Co^{r}(W_j)}\leq C$.
Moreover, usual distortion estimates show that $|e^{S_n
\phi}|_{\Co^{r-1+\coe} (W_j)} \leq C |e^{S_n \phi}|_{\Co^0(W_j)}$,
and that $|e^{S_{n-2L} \phi}|_{\Co^0(T^L W_j)} \leq C \inf_{x\in
B_j} e^{S_{n-2L} \phi(x)}$. Using these estimates, we can compute:
  \begin{align*}
  | e^{S_n \phi} \pi_n|_{\Co^{r-1+\coe}(W_j)}&
  \leq C |e^{S_n \phi}|_{\Co^0(W_j)}
  \leq C |e^{S_{n-2L} \phi}|_{\Co^0(T^L W_j)}
  \\&
  \leq C |e^{S_{n-2L} \phi}|_{\Co^0(T^L W_j)} \int_{B_j} \alpha
  \leq C \int_{B_j} e^{S_{n-2L} \phi} \alpha
  \\&
  = C \int_{B_j} e^{S_{n-2L} \phi} \trunc_{n-2L} \alpha
  = C \int_{T^{n-2L}B_j} \Lp_{\trunc,\phi}^{n-2L} \alpha.
  \end{align*}
Summing over $j$ and using the fact that there is a bounded number
of overlap,
  \begin{equation}
  \sum_j  | e^{S_n \phi} \pi_n|_{\Co^{r-1+\coe}(W_j)}
  \leq C \int_{O\cap T^{-L} W^e} \Lp_{\trunc,\phi}^{n-2L} \alpha.
  \end{equation}
Since the set of integration can be covered by a uniformly bounded
number of admissible leaves, we get
  \begin{equation}
  \sum_j  | e^{S_n \phi} \pi_n|_{\Co^{r-1+\coe}(W_j)}
  \leq C \norm{\Lp_{\trunc,\phi}^{n-2L}\alpha}_{p,q,\coe}.
  \qedhere
  \end{equation}
\end{proof}

\begin{cor}
\label{cor:lowerbound} Assume that $(p,q,\coe)$ is correct. The
spectral radius of $\Lp_{\trunc,\phi}$ acting on $\B^{p,q,\coe}$ is
exactly $\ra$.
\end{cor}
\begin{proof}
Choose once and for all an element $\betriem$ of $\Se$ induced by a
Riemannian metric, as explained in Remark \ref{rem:MetricsinB}. It
satisfies the assumptions of Lemma \ref{lem:Boundsra}. Therefore,
for some constants $L>0$ and $C>0$,
  \begin{equation}
  \label{eq:Boundsran}
  \ra_n^n \leq C \norm{ \Lp_{\trunc,\phi}^{n-2L} \betriem}_{p,q,\coe}
  \leq C \norm{ \Lp_{\trunc,\phi}^{n-2L}}_{p,q,\coe}.
  \end{equation}
Letting $n$ tend to infinity, we obtain that the spectral radius of
$\Lp_{\trunc,\phi}$ is at least $\limsup \ra_n =\ra$. The result
follows remembering Proposition \ref{thm:SpectralGap}.
\end{proof}

\subsection{First description of the peripheral eigenvalues}
\label{sec:PerSpectrum}

In this paragraph, we will study the eigenvalues of modulus $\ra$.
The main goal is to prove that the eigenfunctions for eigenvalues of
modulus $\ra$ are in fact measures. Fix a correct $(p,q,\coe)$.

Denote by $(\gamma_i \ra)_{i=1}^M$ the peripheral eigenvalues of
$\Lp_{\trunc,\phi}$ acting on $\B^{p,q,\coe}$, with $|\gamma_i|=1$.
Let $\kappa$ be the size of the largest Jordan block. Since
$\Lp_{\trunc,\phi}:\B^{p,q,\coe}\to\B^{p,q,\coe}$ is quasicompact,
it must have the form
  \begin{equation}
  \Lp_{\trunc,\phi}=\sum_{i=1}^M(\gamma_i\ra
  S_{\gamma_i}+N_{\gamma_i})+R
  \end{equation}
where $S_{\gamma_i},N_{\gamma_i}$ are finite rank operators such
that $S_{\gamma_i}S_{\gamma_j}=\delta_{ij}S_{\gamma_i}$,
$S_{\gamma_i}N_{\gamma_j}=N_{\gamma_j}S_{\gamma_i}=\delta_{ij}N_{\gamma_j}$,
$N_{\gamma_i}N_{\gamma_j}=\delta_{ij}N_{\gamma_i}^2$,
$S_{\gamma_i}R=RS_{\gamma_i}=N_{\gamma_i}R=RN_{\gamma_i}=0$,
$N_{\gamma_i}^\kappa=0$, and $R$ has spectral radius strictly
smaller than $\ra$. Accordingly, for each $|\gamma|=1$, holds
  \begin{equation}\label{eq:beta}
  \lim_{n\to\infty}
  n^{-\kappa}\sum_{k=0}^{n-1}\gamma^{-k}\ra^{-k}\Lp_{\trunc,\phi}^k
  =\frac{1}{\kappa!} \sum_{i=1}^{M}N_{\gamma_i}^{\kappa-1}\delta_{\gamma
  \gamma_i}.
  \end{equation}
In this formula, if $\kappa=1$, then $N_{\gamma_i}^{\kappa-1}$
indicates the eigenprojection corresponding to the eigenvalue
$\gamma_i \ra$, i.e., $S_{\gamma_i}$. We will denote by
$F_{\gamma_i}$ the image of $N_{\gamma_i}^{\kappa-1}$.

\begin{lem}
\label{lem:rakappaOK} There exists $C>0$ such that, for all $n>0$,
  \begin{equation}
  \ra_n^n \leq C n^{\kappa-1} \ra^n.
  \end{equation}
\end{lem}
This lemma implies in particular that we can apply Lemma
\ref{lem:LY} to $(\ora,d)=(\rho,\kappa-1)$.
\begin{proof}
There exists a constant $C>0$ such that
$\norm{\Lp_{\trunc,\phi}^n}_{p,q,\coe} \leq C n^{\kappa-1} \ra^n$.
Equation \eqref{eq:Boundsran} then implies $\ra_n^n \leq C
n^{\kappa-1} \ra^n$.
\end{proof}

\begin{lem}
\label{lem:AreMeasures0} For all $\gamma$ with $|\gamma|=1$, and all
$\alpha \in F_\gamma$, there exists $C>0$ such that, for all $W\in
\Sigma$, for all $t\leq p$, for all $v_1,\dots, v_t\in
\vectfield^{q+t-\coe}(W)$ with $|v_i|_{\Co^{q+t-\coe}}\leq 1$, for
all $\vf \in \Co^{q+t}_0(W)$,
  \begin{equation}
  \left| \int_W \vf \cdot L_{v_1}\dots L_{v_t} \alpha \right|
  \leq C |\vf|_{\Co^t(W)}.
  \end{equation}
\end{lem}
The point of this lemma is that the upper bound depends only on
$|\vf|_{\Co^t}$ while the naive upper bound would use
$|\vf|_{\Co^{q+t}(W)}$.
\begin{proof}
We can apply Lemma \ref{lem:LY} (and more precisely the inequality
\eqref{eq:bounded}) to $(\ora,d)=(\ra, \kappa-1)$, and to the
parameters $(t,0,\coe)$. We get
  \begin{equation}
  \label{eq:Samekappa}
  \norm{\Lp_{\trunc,\phi}^n}_{t,0,\coe} \leq C n^{\kappa-1} \ra^n.
  \end{equation}

Since $\Se$ is dense in $\B^{p,q,\coe}$, we have
$N_{\gamma}^{\kappa-1} \B^{p,q,\coe}=N_{\gamma}^{\kappa-1} \Se$.
Therefore, we can write $\alpha$ as
$N_{\gamma}^{\kappa-1}(\tilde\alpha)$ where $\tilde\alpha\in \Se$.
Then, by \eqref{eq:beta},
  \begin{equation}
  \label{eq:limit}
  \int_W \vf \cdot L_{v_1}\dots L_{v_t} \alpha = \lim_{n\to\infty}
  \frac{\kappa!}{n^\kappa}\sum_{k=0}^{n-1} (\gamma \ra)^{-k} \int_W \vf \cdot
  L_{v_1}\dots L_{v_t}(\Lp^k_{\trunc,\phi} \tilde\alpha).
  \end{equation}
Moreover, these integrals satisfy
  \begin{equation}
  \left| \int_W \vf \cdot
  L_{v_1}\dots L_{v_t}(\Lp^k_{\trunc,\phi} \tilde\alpha) \right|
  \leq |\vf|_{\Co^t(W)} \norm{ \Lp^k_{\trunc,\phi}
  \tilde\alpha}_{t,0,\coe},
  \end{equation}
by definition of $\norm{\cdot }_{t,0,\coe}$ (this last norm is well
defined since $\tilde\alpha \in \Se$). Using the inequality
\eqref{eq:Samekappa} and the last two equations, we get the lemma.
\end{proof}

Choose $\betriem$ as in Corollary \ref{cor:lowerbound} and let
$\betlim:=\frac{1}{\kappa!} N^{\kappa-1}_1\betriem$. Clearly
$\Lp_{\trunc,\phi}\betlim=\ra\betlim$.
\begin{lem}
\label{lem:AreMeasures} Assume that $(p,q,\coe)$ is correct and
$p>0$. Take $\gamma$ with $|\gamma|=1$, and $\alpha \in F_\gamma$,
Then, for each $W\in\Sigma$, $\alpha$ defines a measure on $W$. In
addition, all such measures are absolutely continuous, with bounded
density, with respect to the one induced by $\betlim$.
% moreover $\betlim\neq 0$.
%Finally, setting $\Lambda^+_n:=\bigcap_{k=0}^{n-1}
%T^k\overline{U'}$, for all $n\in\N$ holds
%$\Id_{\Lambda^+_n}\betlim\equiv \betlim$.\footnote{By this, we mean: for
%all $\Co^\infty$ function $\vf$ which is equal to $1$ on
%$\Lambda^+_n$, $\vf\betlim=\betlim$. Letting $n$ tend to $\infty$, we get
%that, for any $\Co^\infty$ function $\vf$ which is equal to $1$ on a
%neighborhood of $\Lambda^+_\infty$, $\vf\betlim=\betlim$.
%That is, $\betlim$ is
%``supported" on the set $\Lambda^+_\infty$.}
\end{lem}
\begin{proof}
For $t=0$, Lemma \ref{lem:AreMeasures0} shows that $\left|\int_W \vf
\alpha\right| \leq C |\vf|_{\Co^0}$. This shows that $\alpha$
induces a measure on each $W\in \Sigma$.

%For some $\tilde\alpha \in \Se$, we have
%  \begin{equation}
%  \int_W \vf\alpha= \lim_{n\to\infty} \frac{\kappa!}{n^\kappa} \sum_{k=0}^{n-1} (\gamma
%  \ra)^{-k} \int_W \vf \Lp_{\trunc,\phi}^k \tilde\alpha.
%  \end{equation}
For $\gamma=1$ and $\tilde\alpha=\betriem$, $t=0$, Equation
\eqref{eq:limit} shows that $\betlim$ is a nonnegative measure.
Moreover, whenever $\vf\in \Co^q(W)$, it also implies
\[
\left| \int_W \vf \alpha \right| \leq C \int_W |\vf| \betlim.
\]
This inequality extends to continuous functions by density. Hence,
the measure defined by $\alpha$ is absolutely continuous with
respect to the one defined by $\betlim$ (with bounded density).
%
%Accordingly, if $\betlim=0$, then there would be no eigenvalues of modulus $\ra$.
%This would contradict Corollary \ref{cor:lowerbound}
%
%
%To verify the last property suppose $\psi\in\Co^0(U')$ and
%$\psi\Id_{T^nU'}=0$, then, for each $W\in\Sigma$ and
%$\vf\in\Co^{r}_0(W)$, holds
%  \begin{equation}
%  \int_W\vf\psi\betlim=\sum_j\int_{W_j}\vf\circ T^n\rho_j\psi\circ
%  T^n\ra^{-n}e^{S_n(\phi)}\trunc_n\betlim.
%  \end{equation}
%The result then follows since $|\psi|\circ T^n\trunc_n\leq
%|\psi|\circ T^n\Id_{U'}=|\psi\Id_{T^n U'}|\circ T^n=0$ and the fact
%that $\betlim$ is a measure on each manifold, hence, by a standard
%approximation argument, $(1-\Id_{T^n U'})\betlim=0$ for each $n\in\N$.
\end{proof}
\label{par:DefinesMbeta}

An element $\alpha$ of $F_\gamma$ defines a measure on each element
of $\Sigma$. Moreover, if $W$ and $W'$ intersect, and $\vf \in
\Co^q$ is supported in their intersection, then $\int_W \vf
\alpha=\int_{W'} \vf \alpha$. Indeed, this is the case for any
element of $\B^{p,q,\coe}$, since it holds trivially for an element
of $\Se$, and $\Se$ is dense in $\B^{p,q,\coe}$. Therefore, the
measures on elements of $\Sigma$ defined by an element of $F_\gamma$
match locally, and can be glued together: if an oriented submanifold
of dimension $d_s$ is covered by elements of $\Sigma$, then an
element of $F_\gamma$ induces a measure on this submanifold. We will
denote by $\M\alpha$ the measure induced by $\alpha$ on each
oriented stable leaf in $U$.

\begin{lem}
\label{lem:Minj} The map $\alpha \mapsto \M\alpha$ is injective on
each set $F_\gamma$. Moreover, $\betlim\not=0$.
\end{lem}
\begin{proof}
Let $\alpha\in F_\gamma$ satisfy $\M \alpha=0$, we will first prove
that
  \begin{equation}
  \label{lem:EquivNormsSuppr}
  \norm{\alpha}_{0,q,\coe}=0.
  \end{equation}

Notice first that Lemma \ref{lem:CloseLeaves} shows that, if $W'\in
\Sigma$ is $(C_\Sigma,\ve)$-close to an element $W$ of $\Sigma$
contained in a stable manifold, then
  \begin{equation}
  \label{eq:CloseLeaves}
  \left| \int_{W'} \vf \alpha\right| \leq C \ve |\vf|_{\Co^q(W')}.
  \end{equation}
Indeed, the assumption $\M \alpha=0$ shows that, for any $\vf \in
\Co^q_0(W)$, $\ell_{(W,\vf)}(\alpha)=0$.

Take now $W\in \Sigma$ and $\vf\in \Co^q_0(W)$. Using the partition
of unity on $T^{-n} W^{(n)}$ given by the definition of admissible
leaves, we get
  \begin{equation}
  \int_W \vf \alpha=\int_W \vf (\gamma \rho)^{-n}
  \Lp_{\trunc,\phi}^n \alpha
  = (\gamma \rho)^{-n}
  \sum_{j=1}^{k} \int_{W_j} \vf\circ T^n \rho_j\pi_n
  e^{S_n \phi} \cdot \alpha.
  \end{equation}
Each $W_j$ is $(C_\Sigma,\ve_n)$-close to an element of $\Sigma$
contained in a stable leaf, where $\ve_n \to 0$ is given by the
definition of admissible sets of leaves. Hence,
\eqref{eq:CloseLeaves} shows that $\left|\int_W \vf\alpha\right|$ is
bounded by
  \begin{equation}
  \ra^{-n} \sum_{j=1}^{k} C \ve_n |\pi_n e^{S_n \phi}|_{\Co^q(W_j)}
  \leq C \ra^{-n} \ra_n^n \ve_n.
  \end{equation}
The sequence $\ra^{-n} \ra_n^n$ grows at most subexponentially,
while $\ve_n$ goes exponentially fast to $0$ by Definition
\ref{DefAdmissibleLeaves}. Therefore, this quantity goes to $0$,
hence \eqref{lem:EquivNormsSuppr}.

%Let us now prove that
%there exists a constant $C>0$ such that, for all $\alpha\in F_\gamma$,
%$\norm{\alpha}_{p,q,\coe}\leq C \norm{\alpha}_{0,q,\coe}$. Together
%with \eqref{lem:EquivNormsSuppr}, this will show that
%$\norm{\alpha}_{p,q,\coe}=0$ and conclude the proof.

Next, if $\alpha\in F_\gamma$, then $\Lp^n_{\trunc, \phi}
\alpha=(\gamma \ra)^n \alpha$. Using the Lasota-Yorke inequality
\eqref{eq:bounded3} (applied to $(\rho,\kappa-1)$ by Lemma
\ref{lem:rakappaOK}), we get for some $\sigma<1$
  \begin{equation}
  \norm{\alpha}_{p,q,\coe}=\ra^{-n} \norm{\Lp_{\trunc,\phi}^n
  \alpha}_{p,q,\coe}
  \leq  C \sigma^n \norm{\alpha}_{p,q,\coe} ,
  \end{equation}
since $\norm{\alpha}_{0,q+p,\coe}\leq \norm{\alpha}_{0,q,\coe}=0$ by
\eqref{lem:EquivNormsSuppr}. Choosing $n$ large yields
$\norm{\alpha}_{p,q,\coe}=0$.

Let us now prove $\betlim \not=0$. Otherwise, $\M \betlim=0$. For
any $\alpha \in F_\gamma$, the measure $\M \alpha$ is absolutely
continuous with respect to $\M\betlim$, hence zero. By injectivity
of the map $\alpha \mapsto \M \alpha$, we get $\alpha=0$. Therefore,
there is no eigenfunction corresponding to an eigenvalue of modulus
$\ra$. This contradicts Corollary \ref{cor:lowerbound}.
\end{proof}

\section{Peripheral Spectrum and Topology}
\label{sec:TopDyn} In this section we establish a connection between
the peripheral spectrum of the operator and the topological
properties of the dynamical systems at hand.

\subsection{Topological description of the dynamics}
\label{subsec:TopDyn}

Let us recall the classical \emph{spectral decomposition} of a map
$T$ as above (see e.g.\ \cite[Theorem 18.3.1]{katok}). Assume that
$T:U\to X$ is a diffeomorphism and that $\Lambda=\bigcap_{n\in\Z}
T^{n}U$ is a compact locally maximal hyperbolic set. Then there
exist disjoint closed sets $\Lambda_1,\dots,\Lambda_m$ and a
permutation $\sigma$ of $\{1,\dots,m\}$ such that $\bigcup_{i=1}^m
\Lambda_i=NW(T_{\upharpoonright_\Lambda})$, the nonwandering set of
the restriction of $T$ to $\Lambda$. Moreover,
$T(\Lambda_i)=\Lambda_{\sigma(i)}$, and when $\sigma^k(i)=i$ then
$T^k_{\upharpoonright_{\Lambda_i}}$ is topologically mixing, and
$\Lambda_i$ is a compact locally maximal hyperbolic set for $T^k$.

Hence, to understand the dynamics of $T$ on $\Lambda$ (and
especially its invariant measures) when
$\Lambda=NW(T_{\upharpoonright_\Lambda})$, it is sufficient to
understand the case when $T_{\upharpoonright_\Lambda}$ is
topologically mixing.

To deal with orientation problems, we will in fact need more than
mixing. Let $\bar\Lambda$ be the set of pairs $(x,E)$ where $x\in
\Lambda$ and $E\in \G$ is $E^s(x)$ with one of its two possible
orientations. Let $\bar T : \bar \Lambda \to \bar \Lambda$ be the
map induced by $DT$ on $\bar \Lambda$, and let $\pr : \bar\Lambda
\to \Lambda$ be the canonical projection. We have a commutative
diagram
  \begin{equation*}
  \xymatrix{
  \bar\Lambda \ar[d]_\pr \ar[r]^{\bar T} &
       \bar\Lambda \ar[d]^\pr \\
  \Lambda \ar[r]^{T} & \Lambda}
  \end{equation*}
Moreover, the fibers of $\pr$ have cardinal exactly $2$. When $T$ is
topologically mixing, there are exactly three possibilities:
\begin{itemize}
\item Either $\bar T$ is also topologically mixing. In this case, we
say that $T$ is \emph{orientation mixing}.
\item Or there is a decomposition $\bar \Lambda=\bar\Lambda_1 \cup
\bar\Lambda_2$ where each $\bar\Lambda_i$ is invariant under $\bar
T$, and the restriction of $\pr$ to each $\bar\Lambda_i$ is an
isomorphism. We say that $T$ is \emph{mixing, but orientation
preserving}.
\item Or there is a decomposition $\bar \Lambda=\bar\Lambda_1 \cup
\bar\Lambda_2$ such that $\pr$ is an isomorphism on each
$\bar\Lambda_i$, and $\bar T$ exchanges $\bar\Lambda_1$ and
$\bar\Lambda_2$. In this case, $T^2$ is orientation preserving as
defined before.
\end{itemize}

To understand the spectral properties of $T$, it is sufficient to
understand the first two cases, since the last one can be reduced to
the second one by considering $T^2$.

In the second case, there exists an orientation of the spaces
$E^s(x)$ for $x\in \Lambda$, which depends continuously on $x$, and
is invariant under $DT$. Let us say arbitrarily that this
orientation is positive. Consequently, if the neighborhood $U$ of
$\Lambda$ is small enough, there exists a decomposition of $\{ (x,E)
\st x\in U, E \in \G \text{ with } E\subset C_s(x)\}$ into two
disjoint sets $S_+$ and $S_-$, the first one corresponding to vector
spaces $E$ whose orientation is close to the positive orientation of
a nearby set $E^s(x)$, and the other one corresponding to the
opposite orientation. The sets $S_+$ and $S_-$ are invariant under
the action of $DT$. Let $\B^{p,q,\coe}_\pm$ denote the closure in
$\B^{p,q,\coe}$ of the elements of $\Se$ which vanish on $S_\mp$.
Then
  \begin{equation*}
  \B^{p,q,\coe} = \B^{p,q,\coe}_+ \oplus \B^{p,q,\coe}_-.
  \end{equation*}
The transfer operator $\Lp_{\trunc,\phi}$ leaves invariant the sets
$\B^{p,q,\coe}_+$ and $\B^{p,q,\coe}_-$. Moreover, there is a
natural isomorphism from $\B^{p,q,\coe}_+$ to $\B^{p,q,\coe}_-$
(corresponding to reversing the orientation), which conjugates the
action of $\Lp_{\trunc,\phi}$ on $\B^{p,q,\coe}_+$ and
$\B^{p,q,\coe}_-$. Hence, the spectral data of $\Lp_{\trunc,\phi}$
acting on $\B^{p,q,\coe}$ are simply twice the corresponding data
for the corresponding action on $\B^{p,q,\coe}_+$. Therefore, when
$T$ is mixing but orientation preserving, we can restrict ourselves
to the study of $\Lp_{\trunc,\phi}$ acting on $\B^{p,q,\coe}_+$.

\subsection{The peripheral spectrum in the topologically mixing case}

\label{subsec:PeripheralTopology} In this paragraph we will assume
that the dynamics has no wandering parts, that is
$NW(T_{\upharpoonright_\Lambda})=\Lambda$. Given the discussion of
the previous section we can thus restrict ourselves to the mixing
case. Under such an assumption we obtain a complete characterization
of the peripheral spectrum. Note that the proof of the next theorem
relies on some general properties of conformal leafwise measures
that, for the reader's convenience, are proved in Section
\ref{sec:LeafMeasures}.
\begin{thm}
\label{thm:PeripheralSpectrum} Assume that $T$ is orientation mixing
(respectively mixing but orientation preserving). Consider the
operator $\Lp_{\trunc,\phi}$ acting on $\B^{p,q,\coe}$ (resp.\
$\B^{p,q,\coe}_+$). Then  $\ra$ is a simple eigenvalue, and there is
no other eigenvalue of modulus $\ra$.
\end{thm}
\begin{proof}
We give the proof e.g.~for the orientation mixing case, the other
one is analogous.

Let us first prove that $\kappa=1$, that is, there is no Jordan
block. We will show that $\betlim$ satisfies the assumptions of
Lemma \ref{lem:Boundsra}. Assume on the contrary that there exists a
small ball $B$ on which the integral of $\betlim$ vanishes, centered
at a point of $\bigcap_{n\in\N} T^n V$. The preimages of such a
small ball accumulate on the stable manifolds of $T$. By invariance,
the integral of $\betlim$ still vanishes on $T^{-n}B$. Taking a
subsequence and passing to the limit, we obtain a small ball $B'$ in
a stable manifold, centered at a point of $\Lambda$, on which
$\betlim=0$. There is a point $x$ in $\Lambda\cap B'$ such that
$\{T^{-n} x\}$ is dense in $\Lambda$. Let $\ve>0$ be such that the
measure $\M \betlim$ induced by $\betlim$ (as defined in Paragraph
\ref{par:DefinesMbeta}) vanishes on $B(x,\ve)$. Using the invariance
of $\betlim$ and the expansion properties of $T^{-n}$, this implies
that $\M \betlim=0$ on each ball $B(T^{-n}x, \ve)$. By continuity
and density, $\M\betlim=0$. This is in contradiction with Lemma
\ref{lem:Minj}.

Therefore, we can apply Lemma \ref{lem:Boundsra} to $\betlim$, and
get $\ra_n^n \leq C \norm{ \Lp_{\trunc,\phi}^{n-2L}
\betlim}_{p,q,\coe}$. Since $\betlim$ is an eigenfunction for the
eigenvalue $\ra$, this yields $\ra_n^n \leq C \ra^n$. The
Lasota-Yorke inequality \eqref{eq:bounded} yields
$\norm{\Lp_{\trunc,\phi}^n}_{p,q,\coe} \leq C \ra^n$. Hence, there
can be no Jordan block.

Let us now prove that $\betlim$ is the only eigenfunction (up to
scalar multiplication) corresponding to an eigenvalue of modulus
$\ra$. Let $\alpha$ be such an eigenfunction, for an eigenvalue
$\gamma \ra$, with $|\gamma|=1$ and $\alpha\not=0$. Notice first
that the leafwise measure $\M\alpha$ is a continuous leafwise
measure, in the sense of Section \ref{sec:LeafMeasures}. Indeed, if
the test function $\vf$ is $\Co^q$, then the continuity property of
leafwise measures is clear for any element of $\Se$, and extends by
density to any element of $\B^{p,q,\coe}$. When $\alpha\in
F_\gamma$, this continuity property extends from $\Co^q$ test
functions to $\Co^0$ test functions by Lemma \ref{lem:AreMeasures}.
Let us check the assumptions of Proposition \ref{prop:EqualMeasures}
(for the map $T^{-1}$). Note first that $T^{-1}$ is topologically
mixing on $\Lambda$ by assumption, and expanding along stable
leaves. Moreover, let $U$ be an open set in a stable leaf,
containing a point $x\in \Lambda$. Since $T^{-1}$ is transitive,
there exists a nearby point $y$ whose orbit under $T^{-1}$ is dense
in $\Lambda$. The point $z=[x,y]= W^s(x) \cap W^u(y)$ belongs to
$\Lambda\cap U$ if $y$ is close enough to $x$, and its orbit under
$T^{-1}$ is also dense in $\Lambda$. Hence, Proposition
\ref{prop:EqualMeasures} applies, and shows that the measure $\M
\alpha$ is proportional to $\M \betlim$. Since $\M\alpha\not=0$ by
Lemma \ref{lem:Minj}, it follows that $\gamma=1$. Moreover, the
equality $\M\alpha= \gamma'\M\betlim$ implies
$\alpha=\gamma'\betlim$, again by Lemma \ref{lem:Minj}.
\end{proof}

In the course of the above proof, we have showed that $\betlim$
gives a positive mass to each ball in a stable manifold, centered at
a point of $\Lambda$. By compactness of $\Lambda$ and the continuity
properties of $\betlim$, this implies the following useful fact:

For any $\delta>0$, there exists $c_\delta>0$ such that, for any
ball $B(x,\delta)$ in the stable manifold of a point $x\in \Lambda$,
  \begin{equation}
  \label{eq:LowerBoundInt}
  \int_{B(x,\delta)} \betlim \geq c_\delta.
  \end{equation}

\begin{rem}
For the case of unilateral subshifts of finite type, or more
generally when the transfer operator acts on spaces of continuous
functions, there is a much simpler argument to exclude the existence
of Jordan blocks (see \cite{keller:markov} or \cite{baladi:decay}),
which goes as follows.

Assume that the spectral radius of $\Lp$ is $\ra$, and that there
exists an eigenfunction $g>0$ corresponding to this eigenvalue.
Then, for any function $f$, there exists $C>0$ such that $|f|\leq C
g$. Therefore, if the size $\kappa$ of the corresponding Jordan
block is $>1$,
  \begin{equation}
  \frac{1}{n^\kappa} \left| \sum_{k=0}^{n-1} \ra^{-k} \Lp^k f\right|
  \leq C \frac{1}{n^\kappa} \sum_{k=0}^{n-1} \ra^{-k} \Lp^k g\to 0.
  \end{equation}
Hence, $\frac{1}{n^\kappa} \sum_{k=0}^{n-1} \ra^{-k} \Lp^k f$
converges to $0$ in the $C^0$ norm. But it converges to the
eigenprojection of $f$ in the strong norm, so this eigenprojection
has to be $0$ for all $f$. This is a contradiction, and $\kappa=1$.

Unfortunately, this simple argument does not apply in our setting
since the elements of our spaces are not functions: even if we have
constructed the analogue of the function $g$, i.e., $\betlim$, there
is no such inequality as $|\alpha| \leq C\betlim$ for a general
$\alpha \in \Se$. This explains why we had to resort to a more
sophisticated proof.
\end{rem}

\section{Invariant measures and the  variational principle}

\subsection{Description of the invariant measure}
\label{sec:description}

In this paragraph, we assume that $T$ is a map on a compact locally
maximal hyperbolic set, which is either orientation mixing, or
mixing but orientation preserving. Choose $p\in \N^*$ and $q>0$ such
that $(p,q,\coe)$ is correct. In the first case, we let $\B=
\B^{p,q,\coe}$ and in the second case $\B=\B^{p,q,\coe}_+$. The
transfer operator $\Lp_{\trunc,\phi}$ acts on $\B$ and has a simple
eigenvalue at $\ra$ and no other eigenvalue of modulus $\ra$, by
Theorem \ref{thm:PeripheralSpectrum}.

Let $\betlim$ be the eigenfunction of $\ra$. The dual operator
acting on $\B'$ also has a simple eigenvalue at $\ra$. Let $\ell_0$
denote the corresponding eigenfunction, normalized so that
$\ell_0(\betlim)=1$.
\begin{lem}
\label{ExistsMeasure} There exists a constant $C>0$ such that, for
all $\vf \in \Co^{r}(U)$, $|\ell_0( \vf \betlim)| \leq C
|\vf|_{\Co^0}$. Moreover, $\ell_0( \vf \betlim)=\ell_0( \vf\circ T
\cdot \betlim)$.
\end{lem}
\begin{proof}
Let us show that, for any $\alpha\in \B$,
  \begin{equation}
  \label{ell_0Bounded}
  |\ell_0(\alpha)| \leq C \norm{\alpha}_{0,p+q,\coe}.
  \end{equation}
Since $\ell_0 = \ra^{-n} {\Lp'}_{\!\!\trunc,\phi}^{\,n} \ell_0$,
  \begin{align*}
  |\ell_0(\alpha)| &= \ra^{-n} |\ell_0( \Lp_{\trunc,\phi}^n \alpha)|
  \leq C \ra^{-n} \norm{ \Lp_{\trunc,\phi}^n \alpha}_{p,q,\coe}
  \\&
  \leq C \ra^{-n} \left[ C \sigma^n \ra^{n}
  \norm{\alpha}_{p,q,\coe} + C\ra^{n}
  \norm{\alpha}_{0,p+q,\coe}\right]
  \end{align*}
for some $\sigma<1$, by \eqref{eq:bounded3}. Letting $n$ tend to
$\infty$, we obtain \eqref{ell_0Bounded}.

Lemma \ref{lem:AreMeasures} for $t=0$ implies that $\norm{\vf
\betlim}_{0,p+q,\coe} \leq C |\vf|_{\Co^0}$. Together with
\eqref{ell_0Bounded}, this leads to $|\ell_0(\vf \betlim)| \leq C
|\vf|_{\Co^0}$.

Finally, we have
  \begin{multline*}
  \ell_0( \vf \betlim) = ( \ra^{-1}\Lp_{\trunc,\vf}' \ell_0)
  (\vf\betlim)
  = \ra^{-1}\ell_0( \Lp_{\trunc,\vf}(\vf \betlim))
  \\
  =\ell_0( \vf\circ T^{-1} \cdot \ra^{-1}\Lp_{\trunc,\phi} \betlim)
  =\ell_0( \vf\circ T^{-1} \cdot \betlim).
  \end{multline*}
This proves the last assertion of the lemma.
\end{proof}

Lemma \ref{ExistsMeasure} shows that the functional
  \begin{equation*}
  \mu : \vf \mapsto \ell_0( \vf \betlim),
  \end{equation*}
initially defined on $\Co^r$ functions, extends to a continuous
functional on continuous functions. Hence, it is given by a
(complex) measure, that we will also denote by $\mu$. Lemma
\ref{ExistsMeasure} also shows that this measure is invariant.
Hence, it is supported on the maximal invariant set in $U$, i.e.,
$\Lambda$.

\begin{lem}
\label{lem:IsPositive} The measure $\mu$ is a (positive) probability
measure.
\end{lem}
\begin{proof}
By equation \eqref{eq:beta}, the subsequent definition of $\betlim$
and Theorem \ref{thm:PeripheralSpectrum} it follows that, for each
$\alpha\in\B^{p,q,\coe}$,
  \begin{equation}
  \lim_{n\to\infty} \ra^{-n} \Lp_{\trunc,\phi}^n \alpha
  = \ell_0(\alpha)\betlim
  \end{equation}
with $\ell_0(\betriem)=1$. Hence, for all $\vf_1,\vf_2\geq 0$ and
$W\in\Sigma$ holds
  \begin{equation}
  0\leq \lim_{n\to\infty}
  \int_W\vf_1\ra^{-n}\Lp_{\trunc,\phi}^n(\vf_2\betriem)
  =\ell_0(\vf_2\betriem)\int_W\vf_1\betlim.
  \end{equation}
We know that the measure defined by $\betlim$ is nonnegative, and
nonzero by Lemma \ref{lem:Minj}. Therefore, there exist $W$ and
$\vf_1$ such that $\int_W\vf_1\betlim>0$. We get, for any $\vf_2\geq
0$, $\ell_0(\vf_2\betriem)\geq 0$. If $\vf\geq 0$, we have (since
$\ell_0$ is an eigenfunction of $\Lp_{\trunc,\vf}'$)
  \begin{align*}
  \ell_0 ( \vf \betlim) & = \lim_{n\to\infty} \ell_0( \vf
  \ra^{-n} \Lp_{\trunc,\phi}^n \betriem)
  = \lim_{n\to \infty} \ell_0( \ra^{-n} \Lp_{\trunc,\phi}^n( \vf\circ
  T^n   \betriem))
  \\&
  = \lim_{n\to\infty} \ell_0( \vf\circ T^n \betriem)
  \geq 0.
  \end{align*}
Hence, the measure $\mu$ is positive. The normalization
$\ell_0(\betlim)=1$ ensures that it is a probability measure.
\end{proof}

Using the spectral information on $\Lp$, we can now prove the
characterization of the correlations for the measure $\mu$ stated in
Theorem \ref{DescribesCorrelations}. This concludes the proof of
Theorem \ref{DescribesCorrelations} provided one shows that $\mu$ is
indeed the unique Gibbs measure, this will be done in Theorem
\ref{thm:IsGibbs}.
\begin{proof}[Proof of Theorem \ref{DescribesCorrelations}]
We will first describe an abstract setting which implies the
conclusion of the theorem, and then show that hyperbolic maps fit
into this setting.

Let $T$ be a map on a space $X$, preserving a probability measure
$\mu$. Let $\F_1$ and $\F_2$ be two spaces of functions on $X$.
Assume that there exist a Banach space $\B$, a continuous linear
operator $\Lp : \B\to \B$ and two continuous maps $\Phi_1 : \F_1 \to
\B$ and $\Phi_2 : \F_2 \to \B$ such that, for all $n\in \N$, for all
$\psi_1\in \F_1$ and $\psi_2\in \F_2$,
  \begin{equation}
  \label{eq:ToSatisfy}
  \int \psi_1 \cdot \psi_2\circ T^n \dd\mu = \langle \Phi_2(\psi_2),
  \Lp^n \Phi_1(\psi_1) \rangle.
  \end{equation}
Then, for any $\sigma$ strictly larger than the essential spectral
radius of $\Lp$, there exist a finite dimensional space $F$, a
linear map $M$ on $F$, and two continuous maps $\tau_1 : \F_1 \to F$
and $\tau_2 : \F_2 \to F'$ such that \eqref{eq:Correlations} holds.
This is indeed a direct consequence of the spectral decomposition of
the operator $\Lp$.

In our specific setting, we take for $\B$ the Banach space defined
above, $\Lp=\ra^{-1}\Lp_{\trunc,\phi}$, $\F_1$ is the closure of the
set of $\Co^r$ functions in $\Co^p(U)$ and $\F_2$ is the closure of
the set of $\Co^r$ functions in $\Co^q(U)$. On the set of $\Co^r$
functions, define $\Phi_1(\psi_1)=\psi_1 \betlim$, and
$\Phi_2(\psi_2)=\psi_2 \ell_0$. By construction,
\eqref{eq:ToSatisfy} holds. We have to check that $\Phi_1$ and
$\Phi_2$ can be continuously extended respectively to $\F_1$ and
$\F_2$. Let us first prove
  \begin{equation}
  \label{goodnorm}
  \norm{ \psi \betlim}_{p,q,\coe} \leq C |\psi|_{\Co^p}.
  \end{equation}
This will imply that $\Phi_1$ can be extended by continuity to
$\F_1$.

To check \eqref{goodnorm}, consider $t\leq p$, let $W\in \Sigma$,
let $v_1,\dots, v_t \in \vectfield^{q+t-\coe}(W)$ and let $\vf \in
\Co_0^{q+t}(W)$. Then
  \begin{equation}
  \int_W \vf \cdot L_{v_1}\dots L_{v_t}( \psi \betlim)
  =\sum_{A \subset \{1,\dots, t\}} \int_W \vf \left(\prod_{i\in A} L_{v_i}
  \right)\psi \cdot \left( \prod_{i\not\in A}L_{v_i} \right) \betlim.
  \end{equation}
Using Lemma \ref{lem:AreMeasures0} to bound each of these integrals,
we get an upper bound of the form $C |\psi|_{\Co^p}$. This proves
\eqref{goodnorm}.

Let us now extend $\Phi_2$. By \eqref{ell_0Bounded}, for any $\alpha
\in \B$,
  \begin{equation*}
  | \Phi_2(\psi) (\alpha)| \leq C \norm{ \psi \alpha}_{0,p+q,\coe}
  \leq C \norm{\psi \alpha}_{0,q,\coe}
  \leq C |\psi|_{\Co^q} \norm{\alpha}_{0,q,\coe}
  \leq C |\psi|_{\Co^q} \norm{\alpha}_{p,q,\coe}.
  \end{equation*}
Hence, $\norm{\Phi_2(\psi)} \leq C |\psi|_{\Co^q}$. In particular,
$\Phi_2$ can be continuously extended to $\F_2$.

The proof is almost complete, there is just a technical subtlety to
deal with. Since $p$ is an integer, $\F_1=\Co^p(U)$. However, when
$q$ is not an integer, $\Co^r(U)$ is not dense in $\Co^q(U)$, hence
$\F_2$ is strictly included in $\Co^q(U)$. To bypass this technical
problem, we rather use $q'<q$ close enough to $q$ so that
$\sigma>\max( \lambda^{-p}, \nu^{q'})$ (where $\sigma$ is the
precision up to which we want a description of the correlations, as
in the statement of the theorem). Let $\F_2$ be the closure of
$\Co^r(U)$ in $\Co^{q'}(U)$. For $\psi_1\in \F_1$ and $\psi_2 \in
\F_2$, we get as above a description of the correlations, with an
error term at most $C \sigma^n |\psi_1|_{\Co^p(U)}
|\psi_2|_{\Co^{q'}(U)}$. Since $\F_2$ contains $\Co^q(U)$, and
$|\psi_2|_{\Co^{q'}(U)} \leq |\psi_2|_{\Co^q(U)}$, this gives the
required upper bound for all functions of $\Co^q(U)$.
\end{proof}

\subsection{Variational principle}
\label{sec:variational}

We will denote by $B_n(x,\ve)$ the dynamical ball of length $n$ for
$T^{-1}$, i.e.,
  \begin{equation*}
  B_n(x,\ve)=\{y\in U \st \forall\; 0\leq i \leq n-1,\; d(T^{-i}y, T^{-i}x)
  \leq \ve\}.
  \end{equation*}
%By $B_{-n}(x,\ve)$, $n\in \N$, we denote then the dynamical ball for
%$T$, clearly $B_n(T^{n}x,\ve)=T^{n}B_{-n}(x,\ve)$.

\begin{prop}
\label{lem:MeasureDynBalls} For all small enough $\ve >0$, there
exist constants $A_\ve, a_\ve>0$ such that, for all $n\in \N$ and
all $x\in \Lambda$,
  \begin{equation}
 a_\ve e^{S_n \bar\phi(T^{-n}x)}\ra^{-n}\leq \mu(B_n(x,\ve))\leq
\mu(\overline{B_n(x,\ve)})  \leq A_\ve e^{S_n
\bar\phi(T^{-n}x)}\ra^{-n},
  \end{equation}
where $\bar\phi$ is defined by $\bar\phi(y)=\phi(y, E^s(y))$.
\end{prop}
\begin{proof}
Let $\vf$ be a nonnegative $\Co^r$ function supported in
$B_n(x,\ve)$, bounded by one, and equal to one on  $B_n(x,\ve/2)$.
We will prove
  \begin{equation}
  a_\ve e^{S_n \bar\phi(T^{-n}x)}\ra^{-n}\leq \mu(\vf) \leq A_\ve
  e^{S_n \bar\phi(T^{-n}x)}\ra^{-n},
  \end{equation}
which will conclude the proof.

Let $W\in \Sigma$, and let $\vf_0\in\Co_0^{q}(W)$ with
$|\vf_0|_{\Co^{q}(W)} \leq 1$. Then
  \begin{equation}
  \label{eq:smooth-ball}
  \begin{split}
  \int_W \vf_0 \vf \betlim &= \int_W \vf_0 \vf \ra^{-n}
  \Lp_{\trunc, \phi}^n \betlim
  \\
  &=\sum_j \int_{W_j} \rho_j \vf_0\circ T^n  \vf\circ T^n \ra^{-n}
  e^{S_n \phi} \pi_n\cdot
  \betlim,
\end{split}
  \end{equation}
where $\rho_j$ is the partition of unity on $T^{-n}W^{(n)}$ given by
the definition of admissible leaves. Since $\vf$ is supported in
$B_n(x,\ve)$, the number of leaves $W_j$ on which $\vf \circ T^n$ is
nonzero is uniformly bounded. On each of these leaves, $e^{S_n
\phi}$ is bounded by $C e^{S_n \bar \phi(T^{-n}x)}$. It follows
% from Lemma \ref{lem:AreMeasures}
that
  \begin{equation*}
  \left| \int_W \vf_0 \vf \betlim \right| \leq C \ra^{-n} e^{S_n
  \tilde \phi(T^{-n}x)}.
  \end{equation*}
Since this estimate is uniform in $W$ and $\vf_0$, the upper bound
is proven.

For the (trickier) lower bound, we proceed in four steps.

\emph{First step. Let us show that, for any piece $W$ of stable leaf
containing a point $y$ with $d(x,y)<\ve/10$ and $\dist(y, \partial
W) \geq 10\ve$, we have
  \begin{equation}
  \label{eq:LowerBoundW}
  \int_W \vf \betlim \geq C_\ve \ra^{-n} e^{S_n
  \bar\phi(T^{-n} x)}.
  \end{equation}}

Indeed, $T^{-n}W$ contains a disk $D$ centered at a point of
$\Lambda$, of radius $\ve/10$, and contained in $T^{-n}
B_n(x,\ve/2)$. The integral of $\betlim$ on such a disk is uniformly
bounded from below by a constant $C_\ve$ (by
\eqref{eq:LowerBoundInt}), and $\vf \circ T^n=1$ on $D$. Therefore,
  \begin{equation*}
  \int_W \vf\betlim = \int_W \vf \ra^{-n} \Lp_{\trunc,\phi}^n \betlim
  = \ra^{-n} \int_{T^{-n}W} \vf\circ T^n e^{S_n \phi} \pi_n \betlim
  \geq \ra^{-n} \int_D e^{S_n \phi} \pi_n \betlim.
  \end{equation*}
Moreover, $\pi_n\betlim=\betlim$ on $D$ by \eqref{eq:support}, and
$e^{S_n \phi} \geq C e^{S_n \bar\phi(T^{-n}x)}$ on $D$. This proves
\eqref{eq:LowerBoundW}.

\emph{Second step. Let us show that, for any $\delta>0$, there
exists $M=M(\ve,\delta)$ such that, for any $m\geq M$, there exists
$C=C(\ve,\delta,m)$ such that, for any piece $W$ of stable manifold
containing a point $y\in \Lambda $ with $\dist(y,\partial W) \geq
\delta$,
  \begin{equation}
  \int_{T^{-m}W} \vf\betlim \geq C \ra^{-n} e^{S_n
  \bar\phi(T^{-n} x)}.
  \end{equation}}

This is a direct consequence of the topological mixing of $T$ on
$\Lambda$: if $m$ is large enough, then $T^{-m}W$ will contain a
subset $W'$ satisfying the assumptions of the first step. Therefore,
\eqref{eq:LowerBoundW} implies the conclusion.

\emph{Third step. Let $W\in \Sigma$ be a piece of stable manifold
containing a point of $\Lambda$ in its interior. Denote by $W^e$ its
enlargement, as in Definition \ref{DefAdmissibleLeaves}. There
exists $C=C(\ve,W)>0$ such that, for any large enough $p\in \N$,
  \begin{equation}
  \label{eq:ConvergenceOnWe}
  \int_{W^e} \ra^{-p} \Lp_{\trunc,\phi}^p(\vf\betlim)
  \geq C \ra^{-n} e^{S_n\bar\phi(T^{-n}x)}.
  \end{equation}}

To prove this, consider $\{W_j\}$ a covering of $T^{-p} W^{(p)}$ as
in the definition of admissible leaves, and $\rho_j$ the
corresponding partition of unity.

As in the proof of Lemma \ref{lem:Boundsra}, there exists an integer
$L$ with the following property: to each $W_j$, we can associate a
small ball $B(y_j,\delta)$ contained in $T^{-(p-L)} W^e$, at a
bounded distance from $T^L W_j$, with $y_j\in \Lambda$. Increasing
$L$ if necessary (this process does not decrease $\delta$), we can
assume $L\geq M(\ve,\delta)$. Since the balls $B_j$ have a bounded
number of overlaps,
  \begin{equation}
  \int_{W^e} \Lp_{\trunc,\phi}^p(\vf\betlim)
  \geq C \sum_j \int_{B_j} \trunc_{p-L} e^{S_{p-L} \phi}
  \Lp_{\trunc,\phi}^L(\vf\betlim).
  \end{equation}
The function $\trunc_{p-L}$ is equal to $1$ on a neighborhood of the
support of $\betlim$, by \eqref{eq:support}, so we can disregard it.
Moreover, $\inf_{B_j} e^{S_{p-L}\phi} \geq C e^{S_{p-L}(y_j)}$. We
get
  \begin{equation}
  \int_{W^e} \Lp_{\trunc,\phi}^p(\vf\betlim)
  \geq C \sum_j e^{S_{p-L}\bar\phi(y_j)} \int_{T^{-L}B_j} e^{S_L \phi}
  \vf\betlim.
  \end{equation}
The second step applies to each of the sets $B_j$. Since $e^{S_L
\phi}$ is uniformly bounded from below, we obtain
  \begin{align*}
  \ra^n\int_{W^e} \Lp_{\trunc,\phi}^p(\vf\betlim) &
  \geq C e^{S_n\bar\phi(T^{-n}x)}\sum_j e^{S_{p-L}\bar\phi(y_j)}
%  \\&
  \geq C e^{S_n\bar\phi(T^{-n}x)}
  \sum_j \int_{T^L W_j} e^{S_{p-L}\phi} \betlim
  \\&
  \geq C e^{S_n\bar\phi(T^{-n}x)} \int_W
  \Lp_{\trunc,\phi}^{p-L}\betlim
  = C e^{S_n\bar\phi(T^{-n}x)} \int_W \ra^{p-L} \betlim,
  \end{align*}
since $\betlim$ is an eigenfunction of $\Lp_{\trunc,\phi}$.

\emph{Fourth Step. Conclusion.} Fix $W\in \Sigma$ satisfying the
assumptions of the third step. When $p\to \infty$,
$\ra^{-p}\Lp_{\trunc,\phi}^p (\vf\betlim)$ converges to
$\ell_0(\vf\betlim) \betlim=\mu(\vf) \betlim$. Passing to the limit
in \eqref{eq:ConvergenceOnWe}, we obtain
  \begin{equation}
  \mu(\vf) \int_{W^e} \betlim \geq C \ra^{-n} e^{S_n
  \bar\phi(T^{-n}x)}.
  \end{equation}
This is the desired lower bound.
\end{proof}

\begin{thm}
\label{thm:IsGibbs} The spectral radius $\ra$ is equal to the
topological pressure $e^{\Ptop}$ of the function $\bar\phi$. In
addition, the measure $\mu$ is the unique probability measure
satisfying the variational principle
  \begin{equation*}
  h_{\mu}(T)+ \int \bar \phi \dd\mu = \Ptop.
  \end{equation*}
\end{thm}
In other words, $\mu$ is the so-called \emph{Gibbs measure} of
$T:\Lambda \to \Lambda$, corresponding to the potential $\bar \phi$.

\begin{proof}
The theorem is a completely general consequence of Lemma
\ref{lem:MeasureDynBalls}. Indeed, let $T$ be any continuous
transformation on a compact space $\Lambda$ preserving an ergodic
probability measure $\mu$. Let $\bar \phi$ be a function such that
Lemma \ref{lem:MeasureDynBalls} is satisfied, and there exists $C>0$
such that, for any dynamical ball $B=B_n(x,\ve)$, $\sup_B e^{S_n
\bar \phi} \leq C \inf_B e^{S_n \bar \phi}$ (which is satisfied in
our hyperbolic setting since $\bar\phi$ is H\"{o}lder continuous). Then
$\mu$ satisfies the variational principle and is the unique measure
to do so. This result is due to Bowen, and is proved e.g.~in
\cite[Theorem 20.3.7]{katok}. For the convenience of the reader, let
us sketch the proof.

Recall that the definition of the topological pressure of $\bar\phi$
is given by
\[
\Ptop:=\lim_{\ve\to 0}\liminf_{n\to\infty}\frac 1n\ln
S_d(T,\bar\phi,\ve, n) =\lim_{\ve\to 0}\limsup_{n\to\infty}\frac
1n\ln N_d(T,\bar\phi,\ve, n)
\]
where
\[
\begin{split}
&S_d(T,\bar\phi,\ve, n):=\inf\left\{\sum_{x\in E}
e^{S_n\bar\phi(T^{-n}x)}
\st\Lambda\subset\bigcup_{x\in E}B_{n}(x,\ve)\right\}\\
&N_d(T,\bar\phi,\ve, n):=\sup\left\{\sum_{x\in E}
e^{S_n\bar\phi(T^{-n}x)}\st E\subset \Lambda \text{ is }
(n,\ve)\text{-separated}\right\}.
\end{split}
\]

Now in the first case
\[
1=\mu(\Lambda)\leq \sum_{x\in E}\mu(B_n(x,\ve)) \leq
A_\ve\ra^{-n}\sum_{x\in E} e^{S_n\bar\phi(T^{-n}x)}
\]
Taking the inf on $E$ and the limits yields $\ra\leq \Ptop$. On the
other hand if $E$ is $(n,\ve)\text{-separated}$, holds
\[
1=\mu(\Lambda)\geq \sum_{x\in E}\mu(B_n(x,\ve/2)) \geq
a_{\ve/2}\ra^{-n}\sum_{x\in E} e^{S_n\bar\phi(T^{-n}x)}
\]
which, taking the sup on $E$ and the limits, yields $\ra\geq \Ptop$.

Finally, if $\nu$ is any invariant ergodic probability measure, the
Brin-Katok local entropy theorem \cite{brin_katok} states that the
quantity
  \begin{equation*}
  \lim_{\ve \to 0} \limsup_{n\to \infty} \frac{1}{n} \log( 1/
  \nu(B_n(x,\ve)))
  \end{equation*}
converges $\nu$ almost everywhere to $h_\nu(T)$. Lemma
\ref{lem:MeasureDynBalls} shows that, $\mu$-a.e.,
  \begin{equation*}
 \Ptop - \limsup_{n\to\infty} \frac{S_n \bar
  \phi(T^{-n}x)}{n}\geq h_\mu(T)
  \geq \Ptop - \liminf_{n\to\infty} \frac{S_n \bar
  \phi(T^{-n}x)}{n}.
  \end{equation*}
By Birkhoff Theorem, for $\mu$-almost all $x$, $\frac{S_n \bar
  \phi(T^{-n}x)}{n}$ converges to $\int \bar \phi \dd\mu$. Together
with the above inequalities, we get
  \begin{equation*}
  h_\mu(T)+ \int \bar \phi \dd\mu = \Ptop.
  \end{equation*}
Hence $\mu$ maximizes the variational principle. To show that the
maximizing probability is unique one can proceed exactly as in
\cite[Theorem 20.3.7]{katok} where one uses Lemma
\ref{lem:MeasureDynBalls} instead of \cite[Lemma 20.3.4]{katok}.
\end{proof}

\begin{rem}
Theorem \ref{thm:IsGibbs} implies in particular that the measure
$\mu$ constructed using the transfer operator $\Lp_{\trunc,\phi}$ is
in fact independent of the truncation $\trunc$. This can also be
checked directly by spectral arguments. However, $\betlim$ and
$\ell_0$ \emph{do} depend on the truncation: if we take a truncation
with smaller support $\trunc'$, such that $\trunc=1$ on the support
of $\trunc'$, then the new eigenfunctions $\alpha'_0$ and $\ell'_0$
are equal to $\betlim \cdot \prod_{i=1}^{N} \trunc'\circ T^{-i}$ and
$\ell_0 \cdot \prod_{i=0}^{N-1} \trunc'\circ T^i$ for any large
enough $N$. Nevertheless, this shows that they coincide with
$\betlim$ and $\ell_0$ on a neighborhood of $\Lambda$.
\end{rem}

\section{Relationships with the classical theory of Gibbs
measures} \label{sec:realtionships}

\subsubsection{Margulis' construction}
Classically, the Gibbs measure can be constructed by coding, but
there is also a geometric construction, due initially to Margulis.
He proves the following result (for the measure of maximal entropy
in \cite{margulis:these}, but the proofs extend to Gibbs measures,
see e.g.~\cite{babillot_ledrappier}):

There exist a family of measures $\mu^s$ on the stable leaves,
supported on $\Lambda$, and a family of measures $\mu^u$ on unstable
leaves, supported on $\Lambda$, such that
  \begin{equation}
  \label{eq:Conformality}
  \mu^s = T_*(e^{\bar\phi-\Ptop} \mu^s), \quad
  \mu^u=T_*(e^{\Ptop- \bar\phi} \mu^u).
  \end{equation}

The measures $\mu^s$ are constructed by starting from the Riemannian
measure on a very large piece of stable leaf, and then pushing it by
the dynamics $T^n$ (with a suitable multiplication by the weight
$e^{\bar\phi}$). The sequence is shown to converge in some sense, to
the invariant set of measures $\mu^s$. This corresponds exactly to
what we do by the iteration of the transfer operator, exhibiting
$\betlim$ as the limit of $\Lp_{\trunc,\phi}^n (\betriem)$. The main
difference is that we get the convergence in a strong sense (norm
convergence), and for free due to the spectral properties of the
operator. In fact, the measures $\mu^s$ are exactly the measures
induced by $\betlim$ on the stable leaves.

The measures $\mu^u$ are constructed in the same way, but iterating
$T^{-1}$. The relationship with our abstract eigenfunction $\ell_0$
in the dual of $\B$ is less clear at first sight. However, they are
still very closely related. Indeed, let us define an element $\ell
\in \B'$ as follows: if $\alpha \in \B$, and $\vf$ is a $\Co^r$
function supported in a small open set foliated by small stable
leaves, and having as transversal a small unstable leaf $F$, set
  \begin{equation}
  \label{DefineEll}
  \ell( \vf\alpha)=\int_{x\in F} \left( \int_{y\in W^s(x)}
  \vf(y)\prod_{k=0}^\infty \trunc\circ T^k(y)
  e^{\sum_{k=0}^\infty \bar\phi(T^k y)-\bar\phi(T^k x)} \alpha \right)
  \dd\mu^u_F(x).
  \end{equation}
This is well defined since the function $y\mapsto \prod_{k=0}^\infty
\trunc\circ T^k(y) e^{\sum_{k=0}^\infty \bar\phi(T^k y)-\bar\phi(T^k
x)} $ is $\Co^{r-1+\coe}$ on each stable leaf (the product is in
fact finite, since $\trunc\circ T^k$ is uniformly equal to $1$ for
large enough $k$), and can therefore be integrated against $\alpha$.
The Jacobian of the holonomy of the stable foliation with respect to
the measures $\mu^u$ is exactly $e^{\sum_{k=0}^\infty \bar\phi(T^k
y)-\bar\phi(T^k x)}$. Hence, the local definition of $\ell$ is
independent of the choice of the transversal $F$. Using a partition
of unity $\vf_1,\dots,\vf_n$, we have a well defined element
$\ell\in \B'$.
\begin{dft}
Let me recall how this holonomy formula is proved, to be sure that I
have the good signs.

Take $x$ and $y$ on a same stable leaf. Let $A$ and $B$ be small
sets in $W^u(x)$ and $W^u(y)$, corresponding one to each other under
the stable holonomy. The equation $\mu^u= T_*^n (e^{n\Ptop-S_n
\bar\phi} \mu^u)$ shows that $\mu^u(T^n A)$ is close to $e^{n\Ptop
-S_n\bar\phi(x)}\mu^u(A)$, and in the same way $\mu^u(T^n B)$ is
close to $e^{n\Ptop -S_n\bar\phi(y)}\mu^u(B)$. If $n$ is large, $T^n
A$ and $T^n B$ are very close to each other, and have approximately
the same measure. Hence, $\mu^u(B) \sim e^{S_n\bar \phi(y)-S_n \bar
\phi(x)} \mu^u(A)$. We conclude by letting $n$ go to infinity and
$A$ shrink around $x$.
\end{dft}

The conformality property of the measures $\mu^u$ implies that
$\Lp_{\trunc,\phi}' \ell=\ra\ell$. Indeed, let us compute locally:
  \begin{multline*}
  \ell(\Lp_{\trunc,\phi} \alpha)
  =\int_{x\in F} \left( \int_{y\in W^s(x)} \prod_{k=0}^\infty
  \trunc\circ T^k(y) e^{\sum_{k=0}^\infty \bar \phi(T^k y)-\bar
  \phi(T^k x)} \Lp_{\trunc, \phi} \alpha \right) \dd \mu^u_F(x)
  \\
  =\int_{x'\in T^{-1}F} \left( \int_{y'\in W^s(x')} \prod_{k=1}^\infty
  \trunc\circ T^k(y') e^{\sum_{k=1}^\infty \bar \phi(T^k y')-\bar
  \phi(T^k x')} \pi(y') e^{\bar\phi(y')} \alpha \right) \dd
  \mu^u_F(x).
  \end{multline*}
The equality $\mu^u=T_*(e^{\Ptop- \bar\phi} \mu^u)$ gives $\dd
\mu^u_F(x) = e^{\Ptop -\bar \phi(x')} \dd\mu^u_{T^{-1}F}(x')$. It
follows that $\ell(\Lp_{\trunc,\phi} \alpha) = \ra \ell( \alpha)$.

Since the eigenspace of $\Lp'_{\trunc, \phi}$ is one-dimensional,
this shows that $\ell$ and $\ell_0$ are proportional. Hence, the
measures $\mu^u$ give a geometric description of $\ell_0$.

\begin{rem}
This description implies that
  \begin{equation*}
  |\ell_0(\psi \alpha)|\leq C \norm{\alpha}_{0, q,\coe} \cdot
  \sup_{x\in \Lambda} |\psi|_{\Co^{q}(W^s(x))}.
  \end{equation*}
Hence, in \eqref{eq:Correlations}, the factor $|\psi|_{\Co^{q}(U)}$
can be replaced with $\sup_{x\in \Lambda} |\psi|_{\Co^{q}(W^s(x))}$.
\end{rem}

Finally, the Gibbs measure $\mu$ is constructed by ``putting
together locally'' the measures $\mu^s$ and $\mu^u$. In our setting,
this task is automatically performed by the functional analytic
framework.

\subsubsection{Currents}

Another classical construction of Gibbs measures, closely related to
the previous one but expressed slightly differently, is to work with
\emph{currents}, \cite{ruelle_sullivan:currents}. A current of
degree $k$ is an element of the dual of the space of smooth
differential forms of degree $d-k$, where $d$ is the dimension of
the ambient manifold (which we shall assume to be oriented in this
paragraph). A differential form of degree $k$ gives a current of
degree $k$, since it is possible to take its exterior product
against a form of degree $d-k$, and then integrate on the whole
manifold.

A way to construct Gibbs measures is to find ``conformal currents''
in the stable and unstable directions (i.e., currents satisfying a
condition similar to \eqref{eq:Conformality}), and then take their
``intersection'' to get an invariant measure, which is the Gibbs
measure.

Since the differential forms of degree $d_s$ form a subset of $\B$
(see Remark \ref{rem:FormsInB}), an element of the dual of $\B$
gives rise to a current of degree $d_u$. In particular, the
eigenfunction $\ell_0$ is a current (and \eqref{DefineEll} shows
that it is even a current with an interesting underlying geometric
structure). Hence, $\ell_0$ can be interpreted as a conformal
current in the unstable direction.

On the other hand, $\betlim$ is not a current of dimension $d_s$ in
a natural way. Indeed, there is no canonical way to multiply an
element of $\Se$ with a differential form to get something which
could be integrated. However, assume that the weight $\phi$ belongs
to $\Weight^1$ (i.e., it depends only on the point), and that $T$ is
mixing but orientation preserving. Then we can consider in $\B$ the
closure $\CC$ of the set of differential forms. An element of $\CC$
is naturally a current.\footnote{To see this we must check that, if
$\alpha$ is a smooth form of degree $d_u$, there exists $C>0$ such
that, for any form $\beta$ of degree $d_s$, $|\int \alpha \wedge
\beta|\leq C \norm{\beta}_\B$. This can be checked in coordinates by
using a basis of the tangent space whose elements all belong to the
stable cone.} Since $\phi \in \Weight^1$, it is easy to check that
$\Lp_{\trunc, \phi}$ leaves $\CC$ invariant. Moreover, the spectral
radius of the restriction of $\Lp_{\trunc,\phi}$ to $\CC$ is still
$\ra$ (notice that this would \emph{not} hold in the orientation
mixing case). This implies that the eigenfunction $\betlim$ belongs
to $\CC$, hence $\betlim$ can then be interpreted as a current.
Finally, $\mu$ is indeed constructed by ``intersecting'' the two
conformal currents $\ell_0$ and $\betlim$ (this intersection
process, which is often complicated to implement in general, is
given here for free by the functional analytic framework).

\subsubsection{Young-Chernov-Dolgopyat}
In recent years a new approach has been introduced by Lai-Sang
Young. It has been further simplified by Dolgopyat and then
Dolgopyat-Chernov and has been recently reviewed in
\cite{chernov:advanced}. Such an approach is indeed very close to
the one described here. Essentially, it uses objects in the dual of
our spaces $\B^{0,q}$.

More precisely, $\Omega_{p,q,\coe}$, $p+q<r-1+\coe$, can be endowed
with a topology $\tau$, stronger than the weak-* one, for which it
is compact.\footnote{Essentially, two manifolds are close if they
are $\Co^{r-\ve}$ close, for $p+q+\ve<r-1+\coe$, and the $\vf$ must
be $\Co^{p+q-\ve}$ close and the vector fields $\Co^{p+q-\coe-\ve}$
close.} This implies an interesting characterization of the dual
spaces of $\B:=\B^{p,q,\coe}$.
\begin{lem}\label{lem:dual-rap}
Let $\ell_*\in\B'$, then there exists a Borel (with respect to the
$\tau$ topology) measure $\rho$ on $\Omega$ such that, for all
$h\in\B$,
\[
\ell_*(h)=\int_{\Omega}\ell(h)\; \rho(d\ell).
\]
\end{lem}
\begin{proof}
The first step is to construct $F:\B\to\Co^0(\Omega,\C)$ defined by
\[
F(h)(\ell):=\ell(h),
\]
since $\tau$ is stronger than the weak-* topology, $F(h)$ is
continuous. Call $A:=F(\B)$, clearly $A$ is a closed linear space in
$\Co^0(\Omega,\C)$. We can then associate to $\ell_*$ the element
$\nu\in A'$ defined by $\nu(F(h))=\ell_*(h)$. By the Hahn-Banach
Theorem there exists an extension $\nu'$ of $\nu$ to all
$\Co^0(\Omega,\C)$. A this point, by the Riesz representation
Theorem, there exists a measure $\rho$ on $\Omega$ such that
\[
\nu'(f)=\int_{\Omega}f(\ell)\rho(d\ell).
\]
Hence, for each $h\in\B$, we have
\[
\ell_*(h)=\nu(F(h))=\nu'(F(h))=\int_{\Omega}F(h)(\ell)\rho(d\ell)
=\int_{\Omega}\ell(h)\rho(d\ell). \qedhere
\]
\end{proof}
Accordingly, the elements $(W,\vf)\in\Omega_{0,q}$ correspond
exactly to the {\sl standard pairs} in \cite{chernov:advanced} and,
by the above Lemma, the basic objects used in
\cite{chernov:advanced} are precisely the elements of $(\B^{0,q})'$.

The difference lies in the technique used to prove statistical
properties: in \cite{chernov:advanced} is used a probabilistic
coupling technique (instead of the functional analytic one) to prove
statistical properties. Such an approach yields much weaker results
than the present one but it needs much less structure and hence it
is amenable to generalizations in the non-uniformly hyperbolic case.

\subsubsection{Gou\"{e}zel-Liverani}
In \cite{gouezel_liverani}, we introduced an approach to study the
SRB measure of an Anosov map. In many respects, it has the same
flavor as the approach in the present paper, with admissible leaves
and norms obtained in a very similar way. There are however two
important differences between the two papers.
\begin{itemize}
\item
On the technical level, the proof of the Lasota-Yorke inequality
\eqref{eq:bounded2} was more complicated since we had not realized
one could use weighted norms.
\item More conceptually, we had not distinguished between what is
specific to the SRB measure and comes from the Riemannian setting,
and what is completely general. In particular, we considered our
spaces $\B^{p,q}$ as spaces of distributions, by integrating in the
transverse direction with respect to Lebesgue measure. This is very
natural in this case since Lebesgue measure is precisely the
transverse measure $\mu^u$ of Margulis, i.e., the eigenelement
$\ell_0$ in the dual space is already given for free at the
beginning. However, this is really a peculiarity of the SRB measure,
that we had to avoid to treat general Gibbs measures. This explains
why we get spaces of generalized differential forms instead of
spaces of distributions.
\end{itemize}

\section{Examples and Applications}
\label{sec:applications}

In this section we try to give an idea of the breadth of the results
by first discussing some natural examples to which it can be applied
and then illustrating an interesting consequence: perturbation
theory.

\subsection{Examples}

\subsubsection{Anosov and Axiom-A}
Clearly the theory applies to any Anosov or, more generally, Axiom-A
system. In particular, it allows to construct and investigate the
SRB measures and the measures of maximal entropy. In this respect
the present work contains an alternative, self contained,
construction yielding the classical results contained in
\cite{bowen}.\footnote{Notice however that we have an additional
smoothness assumption on the weight.} The relation between the
present approach and other, more classical, ones are discussed in
some detail in Section \ref{sec:realtionships}.

\subsubsection{Open systems} Systems of physical interest are often open,
that is the particles can leave the system. This can happen either
with certainty, once they enter in a given region (holes), or
according to some probability distribution $\pi$ (holes in noisy
systems). The first case cannot be treated in the present setting
since the boundaries of the hole introduce discontinuities in the
system but the latter can be treated provided $\pi$ is smooth. For
example, consider an Anosov system $(X,T)$ and the following
dynamics: a point disappears with probability $\pi(x) dx$ and then,
if it has not disappeared, it is mapped by $T$. In this situation a
typical quantity of physical interest is the escape rate with
respect to Lebesgue, that is the rate at which mass leaks out of the
system. If $\phi$ is the potential corresponding to the SRB measure,
then the transfer operator associated to the above dynamics is
simply $\Lp_{\phi,1-\pi}$ and the escape rate is nothing else than
the logarithm of its leading eigenvalue.

\subsubsection{Billiards with no eclipse conditions} An interesting
concrete system to which the present paper applies is the scattering
by convex obstacle with no-eclipse condition (that is the convex
hull of any two scatterers does not intersect any other scatter).
Although the reflection from an obstacle gives rise to singularities
in the Poincar\'{e} section, nevertheless the no-eclipse condition
implies that only points that will leave the system can experience a
tangent collision (corresponding to a singularity), hence there
exists a neighborhood of the set of the points that keep being
scattered forever in which the dynamics is smooth, hence falls in
our setting. See \cite{knauf-sinai} for a pleasant introduction to
such a subject. In particular, one can obtain sharper information on
the spectrum of the Ruelle operator that are available by the usual
coding techniques used in \cite{morita:eclipse1, stoyanov:eclipse,
morita:eclipse2}.

\subsection{An application: smoothness with respect to parameters}

As already mentioned, the present setting easily allows to discuss
the dependence from parameters of various physically relevant
quantities.

Let us make a simple example to illustrate such a possibility. Let
$(X,T_\lambda)$ be a one parameter family of Anosov maps and let
$\phi_\lambda$ be a one parameter family of potentials. Suppose that
$T_\lambda, \phi_\lambda$ are jointly $\Co^r$ in the variable and
the parameter. By applying the perturbation theory in \cite[Section
8]{gouezel_liverani} it follows that the leading eigenvalue and the
corresponding eigenmeasure are smooth in $\lambda$. If, for example,
we are interested in the measure of maximal entropy
($\phi_\lambda=0$ in view of the variational principle given in
Theorem \ref{thm:IsGibbs}), then it follows that, for any $\ve>0$,
the topological entropy $h_\lambda=P_{\text{top}}(0, T_\lambda)$ is
$\Co^{\lfloor r \rfloor -1-\ve}$ (this is obvious, since this
quantity is constant!) and the measure of maximal entropy
$\mu_\lambda$ is a $\Co^{\lfloor r \rfloor -1-\ve}$ function of
$\lambda$ as a function from $\R$ to $\D_r'$ (that is, if viewed as
a distribution of order $r$).

In fact, the formalism makes it possible to easily compute the
derivatives of the various objects involved. We illustrate this
possibility with the following proposition. Write $T_\lambda$ as
$I_\lambda \circ T_0$ where $I_\lambda$ is the flow from time $0$ to
time $\lambda$ of a $\Co^{r-1}$ time dependent vector field $v_t$.
If $v$ is a smooth vector field, denote by $v^s$ and $v^u$ its
projections on the stable and unstable bundles (they are only H\"{o}lder
continuous vector fields), and by $L_v$ its Lie derivative. If
$\Phi$ is a smooth function on $\G$ such that $\Phi(E)$ is independent of
the orientation of $E$, let $\bar\Phi(x)=\Phi(x,E^s(x))$. The
formula \eqref{DefineEll} for $\ell_0$ shows that, for such a
$\Phi$,
  \begin{equation}
  \label{PhiorbarPhi}
  \ell_0(\Phi \alpha_0)= \mu_0(\bar \Phi).
  \end{equation}
\begin{prop}
 Let $A= \bar\phi'_0-\sum_{n=0}^\infty L_{v_0^s}(\bar \phi_0 \circ
T_0^n)$. Then $h'_0= \mu_0( A)$ and, if $\vf$ is a $\Co^1$ test
function,
  \begin{multline*}
  \left.\frac{d \mu_\lambda(\vf)}{d\lambda}\right|_{\lambda=0}
  = \sum_{k=-\infty}^\infty \mu_0( \vf\circ T_0^k (A-h'_0))\\ +
  \sum_{k=-\infty}^{-1} \mu_0( L_{v_0^u}(\vf\circ T_0^k))
  - \sum_{k=0}^\infty \mu_0( L_{v_0^s}(\vf \circ T_0^k)).
  \end{multline*}
\end{prop}
Notice that the sums in this last equation are clearly finite (the
different terms decay to $0$ exponentially fast). Notice also that,
when the potential $\phi_\lambda$ is constant, we get $h'_0=0$ and,
in the same way, $h'_\lambda=0$. This proves that the topological
entropy is locally constant, without using as usual the structural
stability of the map.
\begin{proof}
Due to \eqref{PhiorbarPhi}, we can omit the bars everywhere and work
only with $\phi_0$.

Let us first prove the following formula. If $W$ is a piece of
stable manifold, $v$ is a smooth vector field on a neighborhood of
$W$ and $\vf \in \Co^1_0(W)$, then
  \begin{equation}
   \label{jashfdlkjsahdfkj}
  \int_W \vf L_v \alpha_0 = -\int_W L_{v^s}\vf \cdot \alpha_0.
  \end{equation}
Notice that $L_{v^s}\vf$ makes sense since $v^s$ is not
differentiated here. To prove this, for large $n$ let $v^{s,n}$ and
$v^{u,n}$ be approximations of $v^s$ and $v^u$ as constructed in
footnote \ref{ft:decompose}. Then
  \begin{equation*}
  \int_W \vf L_{v^{u,n}}\alpha_0 = \ra^{-n}\int_W \vf L_{v^{u,n}}
  (\Lp_0^n \alpha_0)
  = \ra^{-n}
  \int_{T^{-n}W^{(n)}} \vf\circ T_0^n L_{T_0^{*n}v^{u,n}}( \pi_n e^{S_n \phi}
  \alpha_0).
  \end{equation*}
Since $T_0^{*n}v^{u,n}$ has norm at most $C \lambda^{-n}$, this last
integral is bounded by
  \begin{equation}
  C \ra^{-n} \ra_n^n \lambda^{-n} \leq C \lambda^{-n},
  \end{equation}
which tends to $0$ when $n\to \infty$. Hence,
  \begin{equation}
  \int_W \vf L_v \alpha_0 = \int_W \vf L_{v^{u,n}}\alpha_0 - \int_W
  L_{v^{s,n}}\vf \cdot \alpha_0 \to -\int_W L_{v^s}\vf \cdot \alpha_0.
  \end{equation}
This proves \eqref{jashfdlkjsahdfkj}. Together with the formula
\eqref{DefineEll} for the fixed point of the dual operator, we get
for any smooth function $\vf$
  \begin{equation}
  \label{eq:deriv}
  \ell_0( \vf L_v \alpha_0)
  = -\mu_0 (L_{v^s} \vf) -\mu_0\left( \vf\sum_{n=0}^\infty
  L_{v^s}(\phi_0\circ T_0^n)\right).
  \end{equation}

Let $\alpha_\lambda$ be the eigenfunction of the operator
$\Lp_\lambda$ associated to $T_\lambda$ and the potential
$\phi_\lambda$, normalized so that $\ell_0(\alpha_\lambda)=1$. Let
$\ell_\lambda$ be the corresponding eigenfunction of the dual
operator, with $\ell_\lambda(\alpha_\lambda)=1$. The measure
$\mu_\lambda$ is given by $\mu_\lambda(\vf)=\ell_\lambda(\vf
\alpha_\lambda)$. The derivative at $0$ of $\Lp_\lambda \alpha$ is
  \begin{equation}
  \Lp_0'\alpha=L_{v_0}( \Lp_0 \alpha) + \Lp_0( \phi'_0 \alpha).
  \end{equation}
Differentiating the equation $\Lp_\lambda \alpha_\lambda=
e^{h_\lambda}\alpha_\lambda$, we get
  \begin{equation}
  \label{eq:alpha0}
  \alpha'_0 = e^{-h_0}\Lp_0 \alpha'_0 + L_{v_0}\alpha_0 + \phi'_0
  \circ T_0^{-1} \alpha_0 -h'_0 \alpha_0.
  \end{equation}
Applying $\ell_0$ to this equation, we get $h'_0=\mu_0(\phi'_0) +
\ell_0(L_{v_0} \alpha_0)$. By \eqref{eq:deriv} applied to $\vf=1$,
we obtain $h'_0= \mu_0(A)$.

Since $\ell_0(\alpha_\lambda)=1$, we have $\ell_0(\alpha'_0)=0$.
Therefore, $(e^{-h_0} \Lp_0)^n \alpha'_0$ converges to $0$
exponentially fast. We can therefore iterate \eqref{eq:alpha0} and
get
  \begin{equation}
  \alpha'_0=\sum_{k=0}^\infty (e^{-h_0} \Lp_0)^k \bigl[
  L_{v_0}\alpha_0 + (\phi'_0\circ T_0^{-1} -h'_0) \alpha_0 \bigr].
  \end{equation}
We can use this expression to compute $\ell_0(\vf \alpha'_0)$ when
$\vf$ is a smooth function. Let $B=-\sum_{n=0}^\infty
L_{v_0^s}(\phi_0\circ T_0^n)$. Using \eqref{eq:deriv} and
$h'_0=\mu_0(A)$, we obtain
  \begin{multline}
  \label{eq:alpha0def}
  \ell_0(\vf \alpha'_0)
  = \sum_{k=0}^\infty \mu_0( \vf\circ T_0^k (B-\mu_0(B)))
  \\
  -\sum_{k=0}^\infty \mu_0 ( L_{v_0^s}(\vf \circ T_0^k))
  +\sum_{k=1}^\infty \mu_0 ( \vf\circ T_0^k ( \phi'_0 -
  \mu_0(\phi'_0))).
  \end{multline}

For any $\alpha$, we have $\ell_\lambda(\Lp_\lambda \alpha) =
e^{h_\lambda} \ell_\lambda \alpha$. Differentiating, we get
  \begin{equation}
  \label{eq:ell0}
  \ell'_0( \alpha)= \ell'_0(e^{-h_0} \Lp_0 \alpha) + \ell_0( L_{v_0}
  e^{-h_0}\Lp_0 \alpha) + \ell_0 ( (\phi_0'-h'_0)\alpha).
  \end{equation}
Since $\ell_\lambda(\alpha_\lambda)=1$, we have
$\ell_0'(\alpha_0)=-\ell_0(\alpha_0')=0$. Therefore, for any
$\alpha$, $\ell'_0( (e^{-h_0}\Lp_0)^k \alpha)$ converges
exponentially fast to $0$. Iterating \eqref{eq:ell0}, we thus get
  \begin{equation}
  \ell'_0( \alpha)= \sum_{k=0}^\infty \ell_0( L_{v_0} (e^{-h_0}
  \Lp_0)^{k+1} \alpha) + \ell_0( (\phi'_0 -h'_0) (e^{-h_0}\Lp_0)^k
  \alpha).
  \end{equation}
Applying this equation to $\alpha = \vf\alpha_0$ where $\vf$ is a
smooth function, and using $L_{v_0} \vf = L_{v_0^s}\vf +
L_{v^u_0}\vf$ as well as \eqref{eq:deriv}, we get
  \begin{multline}
  \label{eq:ell0def}
  \ell'_0( \vf\alpha_0)
  = \sum_{k=-\infty}^0 \mu_0 ( \vf\circ T_0^k (\phi'_0 -
  \mu_0(\phi'_0)))
  \\
  + \sum_{k=-\infty}^{-1} \mu_0( L_{v_0^u}(\vf\circ T_0^k))
  + \sum_{k=-\infty}^{-1} \mu_0 (\vf\circ T_0^k (B-\mu_0(B))).
  \end{multline}
The derivative at $0$ of $\mu_\lambda(\vf)= \ell_\lambda(\vf
\alpha_\lambda)$ is given by $\ell'_0(\vf\alpha_0)+\ell_0(\vf
\alpha'_0)$. Adding \eqref{eq:ell0def} and \eqref{eq:alpha0def}, we
obtain the conclusion of the proposition.
\end{proof}

Other quantities that can be shown to depend smoothly from
parameters are: the rate of decay of correlations and the associated
distributions $\tau_i$ (see Theorem \ref{DescribesCorrelations}),
the variance in the central limit theorem for smooth observables,
the rate function in the large deviation for  observables (at least
in the $\Co^\infty$ case), etc.

\section{Conformal leafwise measures}

\label{sec:LeafMeasures}

%To study the peripheral spectrum of $\Lp_{\trunc,\phi}$, we will work
%with the measures $\M \alpha$ for $\alpha \in F_\gamma$. These
%measures satisfy conformality properties, and we would like to use
%these properties to show that such measures have to be
%proportional. In this section, we develop abstract tools which give
%such results.

This section is formally independent from the rest of the paper, but
it is of course written with the hyperbolic setting in mind.

Let $X$ be a locally compact space, endowed with a $d$-dimensional
lamination structure: there exists an atlas $\{(U,\chart_U)\}$ where
$U$ is an open subset of $X$ and $\chart_U$ is an homeomorphism from
$U$ to a set $D\times K_U$ where $D$ is the unit disk in $\R^d$ and
$K_U$ is a locally compact space. Moreover, the changes of charts
send leaves to leaves, i.e., $\chart_U\circ \chart_V^{-1}( x, y)=
(f(x,y), g(y))$ where defined.

A \emph{continuous leafwise measure} $\mu$ is a family of Radon
measures on each leaf such that, for all chart $(U,\chart_U)$ as
above and all continuous function $\vf$ supported in $U$,
$\int_{\chart_U^{-1}(D\times\{y\})} \vf \dd\mu$ depends continuously
on $y\in K_U$.

Assume that, on each leaf of the lamination, a distance is given,
which varies continuously with the leaf (in the sense that, for any
chart $(U,\chart_U)$ as above, the map from $D\times D \times K_u$
to $\R$ given by $(x,x',y)\mapsto d(\chart_U^{-1}(x,y),
\chart_U^{-1}(x',y))$ is continuous). Consider then an open subset
$Y$ of $X$, with compact closure, and a continuous map $T:Y\to X$
which sends leaves to leaves and expands uniformly the distance:
there exist $\expansion>1$ and $\delta_0>0$ such that, whenever
$x,y$ are in the same leaf and satisfy $d(x,y)\leq \delta_0$, then
$d(Tx,Ty)\geq \expansion d(x,y)$ (in particular, the restriction of
$T$ to $B(x,\delta_0)$ is a homeomorphism). Assume that
$\Lambda:=\bigcap_{n\geq 0}T^{-n} X$ is a compact subset of $X$.

If $x\in Y$, then $T$ is a homeomorphism on a small ball around $x$
in the leaf containing $x$. Hence, it is possible to define the
pullback $T^* \mu$ of any continuous leafwise measure $\mu$. Our
first result is:

\begin{thm}
\label{EqualMeasures} Let $\mu$  be a nonnegative continuous
leafwise measure, and $\nu$ a complex continuous leafwise measure.
Assume that there exists a constant $C>0$ such that, on each leaf,
$|\nu| \leq C \mu$. Moreover, assume that there exists a continuous
function $\trunc$, supported in $Y\cap T^{-1}Y$, positive on
$\Lambda$,  H\"{o}lder continuous on each leaf, such that $\mu= \trunc
T^*\mu$ and $\nu=\gamma\trunc T^*\nu$ for some $\gamma\in \C$ with
$|\gamma|=1$.

Then there exist $c\in \C$ and an open subset $U$ of a leaf,
containing a point of $\Lambda$, such that $\nu=c\mu$ on $U$.
\end{thm}

The proof is essentially a density point argument: there is a small
subset where $\nu$ is very close to a multiple of $\mu$, and pushing
this estimate by $T^N$ for large $N$ we will obtain the result.
Technically, the existence of convenient density points will be
proved using the martingale convergence theorem. Hence, we will
first need to construct good partitions.

Notice first that
  \begin{equation}
  \label{eq:support}
  \text{the leafwise measure $\mu$ is supported on
  $\Lambda$.}
  \end{equation}
Indeed, if a compact set $V$ of a leaf does not intersect $\Lambda$,
then it can be covered by a finite number of open subsets which are
sent in $X\backslash Y$ by some iterate of $T$. The equation
$\mu=\trunc T^*\mu$ then shows that $\mu$ gives zero mass to each of
these open sets.

By compactness of $\Lambda$, there exist $\delta\in (0,\delta_0)$
and $\ve_0>0$ such that, for any $x\in \Lambda$, the ball
$B(x,\delta)$ (in the leaf containing $x$) is contained in
$\{\pi>\ve_0\}$. We fix such a $\delta$ until the end of the proof.

We will say that a subset $A$ of a leaf is \emph{good} if it is open
with compact closure and $\mu(\partial A)=0$.
\begin{lem}
\label{sous_lemme} Let $A$ be a good subset of a leaf, and let
$\ve>0$. There exist good subsets $B$ and $(F_i)_{1\leq i \leq K}$
forming a partition of a full measure subset of $A$, with
$\diam(F_i)\leq \ve$, such that $\mu(B)\leq \mu(A)/2$ and, for all
$i$, there exist $n\in \N$ and $x\in \Lambda$ such that
$B(x,\delta/5) \subset T^n F_i \subset B(x,\delta)$.
\end{lem}
\begin{proof}
Since $\mu(\partial A)=0$, there exists $\eta>0$ such that $V=\{x\in
A, d(x,\partial A)\geq \eta\}$ satisfies $\mu(V)\geq \mu(A)/2$.
Choose $N>0$ such that $\expansion^N \ve>\delta$ and $\expansion^N
\eta>\delta$.

Define a distance $d_N$ on $A$ by $d_N(x,y)=\sup_{0\leq i\leq
N}d(T^i x, T^i y)$. Let $B_N(x,r)$ denote the ball of center $x$ and
radius $r$ for the distance $d_N$. Choose a maximal
$\delta/2$-separated set for the distance $d_N$ in $\Lambda \cap V$,
say $x_1,\dots, x_k$. Then the balls $B_N(x_i, \delta/4)$ are
disjoint, and $T^N(B_N(x_i,\delta/5))=B(T^N x_i, \delta/5)$.
Moreover, $V\cap \Lambda \subset \bigcup B_N(x_i,\delta/2)$.

For each $i$, there exist $a_i\in (\delta/5,\delta/4)$ with
$\mu(\partial B_N(x_i, a_i))=0$, and $b_i\in (\delta/2,\delta)$ with
$\mu(\partial B_N(x_i, b_i))=0$. Define then the sets $F_i$ by
induction on $i$, by
  \begin{equation*}
  F_i= B_N(x_i, b_i) \backslash \left( \bigcup_{j<i} F_j \cup
  \bigcup_{j>i} \overline{B_N(x_i,a_i)} \right).
  \end{equation*}
By construction, the sets $F_i$ are good sets and $B(T^N
x_i,\delta/5) \subset T^N F_i \subset B(T^N x_i, \delta)$. Set
finally $B=A\backslash \bigcup \overline{F}_i$. The sets $F_i$ cover
almost all $V\cap \Lambda$, i.e. almost all $V$ since $\mu$ is
supported on $\Lambda$. This implies that $\mu(B) \leq
\mu(A\backslash V)\leq \mu(A)/2$.
\end{proof}

\begin{lem}
\label{ConstructsGoodPartitions} Let $A$ be a good subset of a leaf,
and let $\ve>0$. There exist good subsets $(F_i)_{i\in \N}$ of $A$,
with $\diam(F_i) \leq \ve$, forming a partition of a full measure
subset of $A$, such that for all $i\in \N$, there exist $n\in \N$
and $x\in \Lambda$ such that $B(x,\delta/5) \subset T^n F_i \subset
B(x,\delta)$.
\end{lem}
\begin{proof}
It is sufficient to apply inductively Lemma \ref{sous_lemme} to $A$,
then $B$, and so on.
\end{proof}

\begin{proof}[Proof of Theorem \ref{EqualMeasures}]

Let us say that a set has ``full $\mu$ measure'' if its intersection
with any leaf has full measure in the usual sense. Let $f=
\frac{\dd\nu}{\dd\mu}$ be the leafwise Radon-Nikodym of $\nu$ with
respect to $\mu$.  It is defined $\mu$ almost everywhere. Since
$|\nu|\leq C \mu$, it satisfies $|f|\leq C$. The equations
$\mu=\trunc T^* \mu$ and $\nu=\gamma\trunc T^*\nu $ show that, for
almost all $x\in \Lambda$, $f(Tx)= \gamma^{-1} f(x)$.

Start from a good set $A$ in a leaf, containing a point of
$\Lambda$. Applying inductively Lemma
\ref{ConstructsGoodPartitions}, we obtain a sequence of finer and
finer partitions $\FF_n$ of a full measure subset of $A$, such that,
for all $F\in \FF_n$, there exists $i\in \N$ and $x\in \Lambda$ such
that $B(x,\delta/5) \subset T^i F \subset B(x,\delta)$, and with
$\diam F \leq 2^{-n}$.

For $\mu$ almost every $x\in A$, there is a well defined element
$F_n(x)\in \FF_n$ containing $x$. Moreover, the martingale
convergence theorem ensures that, for $\mu$ almost every $x$, for
all $\ve>0$,
  \begin{equation}
  \label{TendsToZero}
  \frac{\mu\{ y\in F_n(x) \st |f(y)-f(x)|>\ve\}}{\mu(F_n(x))} \to
  0 \text{ when }n\to\infty.
  \end{equation}
Fix such a point $x$. Let $x_n\in \Lambda$ and $i(n)\in \N$ be such
that $B(x_n,\delta/5) \subset T^{i(n)} F_n(x) \subset
B(x_n,\delta)$. Since $\trunc$ is H\"{o}lder continuous and $\trunc\geq
\ve_0$ on the iterates $T^j F_n(x)$ for all $0\leq j< i(n)$, there
exists a constant $C$ such that, for all $y,z\in F_n(x)$,
  \begin{equation*}
  \prod_{j=0}^{i(n)-1} \trunc(T^j y) \leq C \prod_{j=0}^{i(n)-1} \trunc(T^j
z).
  \end{equation*}
Together with \eqref{TendsToZero} and the equation $\mu=\trunc
T^*\mu$, this gives
  \begin{equation*}
  \frac{\mu\{ y\in T^{i(n)}F_n(x) \st |f(T^{-i(n)}y)-f(x)|>\ve\}}
  {\mu(T^{i(n)}F_n(x))} \to
  0.
  \end{equation*}
Moreover, $f(T^{-i(n)}y)=\gamma^{i(n)}f(y)$, and $\mu(T^{i(n)}
F_n(x)) \leq \mu(B(x_n, \delta))$ is uniformly bounded. Hence, for
all $\ve>0$,
  \begin{equation*}
  \mu\{ y\in T^{i(n)}F_n(x) \st |f(y)-\gamma^{-i(n)} f(x)|>\ve\}
  \to 0.
  \end{equation*}
Since $T^{i(n)}F_n(x)$ contains the ball $B(x_n,\delta/5)$, we get
in particular
  \begin{equation}
  \label{Converges}
  \mu\{ y\in B(x_n,\delta/5) \st |f(y)-\gamma^{-i(n)} f(x)|>\ve\}
  \to 0.
  \end{equation}
Taking a subsequence if necessary, we can assume that $x_n$
converges to a point $x'$ and $\gamma^{-i(n)}$ converges to
$\gamma'\in \C$ with $|\gamma'|=1$. Let $\vf$ be a continuous
function supported in $B(x',\delta/10)$. Extend it to a continuous
function with compact support on nearby leaves. Then
\eqref{Converges} and the inequality $|f|\leq C$ show that
  \begin{equation*}
  \int_{B(x_n, \delta/5)} \vf \dd\nu - f(x)\gamma' \int_{B(x_n,
  \delta/5)} \vf\dd\mu \to 0.
  \end{equation*}
By the continuity properties of $\mu$ and $\nu$, this implies that
  \begin{equation*}
  \int_{B(x', \delta/10)} \vf\dd\nu= f(x)\gamma' \int_{B(x',\delta/10)}
  \vf\dd\mu.
  \end{equation*}
Hence, on the ball $B(x',\delta/10)$, we have $\nu=f(x)\gamma' \mu$.
\end{proof}

\begin{prop}
\label{prop:EqualMeasures} Under the assumptions of Theorem
\ref{EqualMeasures}, assume moreover that the map $T$ is
topologically mixing on $\Lambda$, and that any open set $U$ of a
leaf which contains a point of $\Lambda$ also contains a point of
$\Lambda$ whose orbit is dense. Then there exists $c\in \C$ such
that $\nu=c\mu$. In particular, $\gamma=1$ (or $\nu=0$).
\end{prop}
\begin{proof}
Note first that, if there exists an open subset $U$ of a leaf on
which $\nu$ vanishes, then $\nu=0$ and the theorem is trivial.
Indeed, there exists $x\in U\cap \Lambda$ whose positive orbit under
$T$ is dense in $\Lambda$. Let $r\in (0,\delta)$ be such that
$B(x,r)\subset U$. The conformality of $\nu$ and the expansion
properties of $T$ show that, for any $n\in \N$, $\nu$ vanishes on
$B(T^n x, r)$. Since $\nu$ is continuous, it follows that $\nu=0$ on
$\Lambda$. Since $\nu$ is supported on $\Lambda$, $\nu=0$.

Assume now that $\nu$ is nonzero on each set $U$ as before. Since
$|\nu|\leq \mu$, this implies the same property for $\mu$. By
Theorem \ref{EqualMeasures}, there exists an open set $U$ in a leaf,
containing a point of $\Lambda$, and $c\in \C$ such that $\nu=c\mu$
on $U$. As above, consider $x\in U\cap\Lambda$ whose orbit is dense,
and choose $r\in (0,\delta)$ such that $B(x,r)\subset U$. The
conformality of $\nu$ and $\mu$ shows that, on $B(T^n x, r)$, $\nu=c
\gamma^{-n} \mu$. By continuity of the measures, for any $y\in
\Lambda$, there exists $f(y) \in \C$ such that $\nu=f(y) \mu$ on
$B(y,r)$. Moreover, this $f(y)$ is uniquely defined since $\mu$ is
nonzero on any ball $B(y,r)$, it depends continuously on $y\in
\Lambda$, and it is nonzero by assumption on $\nu$. Finally, $f\circ
T= \gamma^{-1} f$.

Since $T$ is topologically mixing, this implies that $f$ is constant
and $\gamma=1$.
\end{proof}

\bibliography{biblio}
\bibliographystyle{alpha}

\end{document}